\documentclass[pdflatex]{sn-jnl}
\usepackage[numbers,sort&compress]{natbib}

\usepackage{hyperref}
\hypersetup{
    unicode      = {false},
    pdftoolbar   = {true},
    pdfmenubar   = {true},
    pdffitwindow = {true},
    pdfauthor    = {Audet et al.},
    pdftitle     = {bocatmads},
    pdfsubject   = {bocatmads},
    pdfkeywords  = {bocatmads},
    pdfnewwindow = {true},
    colorlinks   = {true},
    linkcolor    = {blue},
    citecolor    = {blue},
    filecolor    = {black},
    urlcolor     = {blue},
    breaklinks   = {true}
} 

\usepackage{booktabs}
\usepackage[utf8]{inputenc}
\usepackage{comment}
\usepackage[english]{babel}
\usepackage[T1]{fontenc}
\usepackage[export]{adjustbox}
\usepackage{booktabs}
\usepackage{multirow}
\usepackage{bm}
\usepackage{graphicx}
\usepackage{epstopdf}
\usepackage{soul}
\usepackage{textcomp}

\usepackage{amsmath}

\DeclareMathOperator*{\argmax}{argmax}

\usepackage{stackengine}
\usepackage{mathtools}
\usepackage{stackrel}
\usepackage{dsfont}

\usepackage{mathrsfs}   
\DeclareFontFamily{U}{rsfs}{\skewchar\font127 }
\DeclareFontShape{U}{rsfs}{m}{n}{<-6> rsfs5 <6-8> rsfs7 <8-> rsfs10}{}

\usepackage{changepage} 
\usepackage{enumitem}

\usepackage{amsfonts}
\usepackage{makecell}
\usepackage{amssymb}
\usepackage{lipsum}
\usepackage{multicol}

\usepackage{silence}
\usepackage[font=small]{caption}

\usepackage{subcaption}
\usepackage{float}
\usepackage[separate-uncertainty=true]{siunitx}
\usepackage{tikz}

\usepackage{forest}
\usepackage{stanli}
\usetikzlibrary{positioning,shapes,shadows,arrows,automata,arrows.meta,decorations, decorations.text,calligraphy}
\usetikzlibrary{hobby, decorations.pathreplacing,calligraphy}
\usetikzlibrary{patterns}
\usetikzlibrary{matrix}
\usepackage{pgfplots}
\pgfplotsset{compat=1.18}
\usetikzlibrary{hobby, decorations.markings, arrows.meta}

\usepackage{pgfplots}
\pgfplotsset{compat=1.18}
\usepgfplotslibrary{groupplots}

 \def\keywordname{{\bfseries Keywords}}

\usepackage{xfrac}
\usepackage{standalone}
\usepackage{yhmath}

\definecolor{Red}{rgb}{1,0,0}
\definecolor{Green}{rgb}{0,.6,0}
\definecolor{Blue}{rgb}{0,0,1}
\definecolor{Black}{rgb}{0,0,0}

\newcommand{\blr}{\color{Black}}
\newcommand{\bl}{\color{Blue}}

\newcommand{\x}{\bm{x}}
\newcommand{\y}{\bm{y}}
\newcommand{\z}{\bm{z}}
\newcommand{\w}{\bm{w}}
\newcommand{\uu}{\bm{u}}
\newcommand{\vv}{\bm{v}}

\newcommand{\quant}{\mathrm{qnt}}
\newcommand{\integer}{\mathrm{int}} 
\newcommand{\continuous}{\mathrm{cont}} 
\newcommand{\cat}{\mathrm{cat}}

\newcommand{\feasible}{\text{\textsc{fea}}}
\newcommand{\infeasible}{\text{\textsc{inf}}}

\newcommand{\A}{\text{\sc a}\xspace}
\newcommand{\B}{\text{\sc b}\xspace}

\newcommand{\Aluminum}{\text{\sc alum}\xspace}
\newcommand{\Steel}{\text{\sc steel}\xspace}
\newcommand{\Composite}{\text{\sc comp}\xspace}
\newcommand{\Wood}{\text{\sc wood}\xspace}

\newcommand{\Square}{\text{\sc square}\xspace}
\newcommand{\Circle}{\text{\sc circle}\xspace}
\newcommand{\Ellipse}{\text{\sc ellipse}\xspace}

\newtheorem{definition}{Definition}

\usepackage{fancyhdr}
\usepackage{csquotes}
\usepackage{xcolor}

\usepackage[ruled,vlined]{algorithm2e}
\DontPrintSemicolon

\usepackage{xspace}
\newcommand{\mads}{{\sf MADS}\xspace}
\newcommand{\catmads}{{\sf CatMADS}\xspace}
\newcommand{\catmadsf}{{\sf CatMADS}\ensuremath{{}^{f}}\xspace}

\newcommand{\catmadsgp}{{\sf CatMADS\textsubscript{GP}}\xspace}

\newcommand{\catmadsgpf}{{\sf CatMADS}\ensuremath{{}_{\mathrm{GP}}^{f}}\xspace}

\newcommand{\mvmads}{{\sf MV-MADS}\xspace}

\newcommand{\nomadbbo}{{\sf NomadBBO}\xspace}
\newcommand{\catsuite}{{\sf Cat-Suite}\xspace}

\usepackage[nameinlink,capitalise,noabbrev]{cleveref}

\title{{\centering Surrogate-based categorical neighborhoods \\ for mixed-variable blackbox optimization
}}

\author{Charles Audet}
\author{Youssef Diouane}
\author{Edward Hall\'e-Hannan\footnotemark}
\author{S\'ebastien Le~Digabel}
\author{Christophe Tribes}
\affil{GERAD and Department of Mathematics and Industrial Engineering, Polytechnique Montr\'eal, Montr\'eal, QC, Canada}
\begin{document}
%
\abstract{
In simulation-based engineering, design choices are often obtained following the optimization of complex blackbox models.
These models frequently involve mixed-variable domains with quantitative and categorical variables.
Unlike quantitative variables, categorical variables lack an inherent structure, which makes them difficult to handle, especially in the presence of constraints.
This work proposes a systematic approach to structure and model categorical variables in constrained mixed-variable blackbox optimization.
Surrogate models of the objective and constraint functions are used to induce problem-specific categorical pseudo-distances.
From these pseudo-distances, surrogate-based neighborhoods are constructed using notions of dominance from bi-objective optimization, jointly accounting for information from both the objective and the constraint functions.
This study addresses the lack of automatic and constraint-aware categorical neighborhood construction in mixed-variable blackbox optimization.
As a proof of concept, these neighborhoods are employed within \catmads, an extension of the \mads algorithm for categorical variables.
The surrogate models are Gaussian processes, and the resulting method is called \catmadsgp.
%
%
{\blr
The method is benchmarked on the \catsuite collection of~$60$ mixed-variable optimization problems and compared with two external mixed-variable solvers.
Data profiles indicate that \catmadsgp achieves the best performance on this collection of problems.
}
}



\maketitle
\begin{NoHyper}
\footnotetext{*Corresponding author
(\nolinkurl{edward.halle-hannan@polymtl.ca}).}
\end{NoHyper}


%



\noindent
{\small
\textbf{Funding:} This research is funded by a Natural Sciences and Engineering Research Council of Canada (NSERC) PhD Excellence Scholarship (PGS D), a Fonds de Recherche du Qu\'ebec (FRQNT) PhD Excellence Scholarship and an Institut de l’\'Energie Trottier (IET) PhD Excellence Scholarship, 
as well as by the NSERC discovery grants  
    RGPIN-2020-04448 (Audet),
    RGPIN-2024-05093 (Diouane)
    and RGPIN-2024-05086 (Le~Digabel).
}

\newpage
\section{Introduction} 
\label{sec:intro}

In most cases, design choices in engineering and machine learning imply the constrained optimization of complex computer simulations that are \textit{expensive-to-evaluate} blackboxes with no analytical expressions.
In addition, these blackboxes may involve different types of variables, such as continuous, integer and categorical variables.
The latter are qualitative and lack structure.
For example, in aerospace engineering, optimal aircraft design involves simulations with choices of materials or propulsion systems~\cite{BuCiDeNaLa2021,SaBaDiLeMoDaNgDe2021}.
Deep learning models are characterized by hyperparameters, such as the activation function or the choice of optimizer, whose values affect the training and validation pipeline~\cite{G-2022-11, G-2024-33, hypernomad_paper}.
Standard optimization methods using gradients and relaxations in the presence {\blr of} blackboxes and categorical variables are not suited for this context.

The motivation of this work is to improve the optimization of categorical variables in a constrained {\blr blackbox setting} where information is limited and constraints may be difficult to handle.
This is done by proposing a systematic approach to structure categorical variables with neighborhoods constructed using surrogate models of the objective and constraint functions.
These neighborhoods can be directly incorporated into mixed-variable methods using categorical neighborhoods.
As a proof-of-concept, they are employed with \catmads~\cite{catmads}, a mixed-variable extension of the Mesh Adaptive Direct Search (\mads) algorithm~\cite{AuDe2006} for constrained blackbox optimization.

\subsection*{Notation}

The notation used follows~\cite{catmads}. Vectors are in bold, and scalars are in normal font.
Superscripts are reserved for types of variables.
Subscripts without parentheses index variables, \textit{e.g.}, $x_1^{\cat}$ is the first categorical variable.
Subscripts with parentheses are reserved for indexing an iteration $k$, \textit{e.g.} $\x_{(k)}$, or for listing points, \textit{e.g.} $\x_{(1)}, \x_{(2)}, \ldots, \x_{(s)}$.

Sets are denoted using capital letters with either normal or calligraphic font, \textit{e.g.} 
A or $\mathcal{A}$, except for the set of known points noted $\mathbb{X}$.

\subsection{Constrained mixed-variable blackbox optimization}
\label{sec:problem_statement}

This study tackles inequality constrained mixed-variable blackbox optimization problems formulated as  
\begin{equation}
    \min_{ \x \in \Omega } f(\x), 
    \label{eq:formulation}
\end{equation}
where $f:\mathcal{X} \to \overline{\mathbb{R}}$ is the objective function with $\overline{\mathbb{R}} = \mathbb{R} \cup \{+ \infty \}$, $\mathcal{X}$ is the domain of the objective and constraint functions, $\x \in \mathcal{X}$ is a point, 
$\Omega \coloneq \{ \x \in \mathcal{X} \ : \ g_j( \x) \leq 0, \  j \in J \}$ is the feasible set defined by the constraints, $g_j:\mathcal{X} \to \overline{\mathbb{R}}$ is the $j$-th constraint function of the problem with $j \in J$ and $|J|\in \mathbb{N}$ is the number of constraints.

The constraints defining the feasible set $\Omega$ are supposed to be \textit{relaxable} and \textit{quantifiable}, meaning that when they are evaluated, they can
return a proper value without being satisfied~\cite{LedWild2015}.
An \textit{unrelaxable} constraint, which must be satisfied to have a proper blackbox execution, is considered out of the domain, such that if $\x \not \in \mathcal{X}$, then it is assigned $f(\x)=+\infty$.
Similarly, if a point $\x \in \mathcal{X}$ hits a \textit{hidden constraint} that crashes or invalidates the output values, then it is simply assigned $f(\x)=+\infty$.

The functions involved in the main optimization problem are supposed to be \textit{blackboxes}, which are defined as follows in~\cite{AuHa2026}: ``\textit{any process that when provided an input, returns an output, but the inner workings of the process are not analytically available}''. 
%
%
%
%
{\blr
For example, in aeronautics, the quality of a design may be evaluated through a complex simulation combining aerodynamic and structural computations~\cite{Saves2024}.
Such a simulation takes aircraft design variables as input, such as its geometry and material choices, and returns performance measures as output.
Derivatives with respect to the continuous variables are often unavailable or even nonexistent.
Moreover, executing such simulations to obtain function evaluations may be computationally expensive.
The evaluation budget is therefore typically limited.
}


The problems considered are said to be mixed-variable.
Each variable is either \textit{continuous} ($\continuous$), \textit{integer} ($\integer$) or \textit{categorical} ($\cat$).
The variables of the same type are contained in a corresponding column vector called a \textit{component}.
For $t \in \{ \cat, \integer, \continuous  \}$, the component of type $t$ is expressed as
\begin{equation*}
    \x^t \coloneq \left(x^t_1, x^t_2, \ldots, x^t_{n^t} \right) \in \mathcal{X}^t \coloneq \mathcal{X}^t_1 \times \mathcal{X}^t_2 \times \ldots \times \mathcal{X}^t_{n^t},
\end{equation*}
where $x_i^t \in \mathcal{X}_i^t$ represents the $i$-th variable of type $t$, and $\mathcal{X}_i^t$ denotes its bounds.
The variable index $i$ ranges over the set $I^t \coloneq \{1,2,\ldots, n^t\}$.
The number of variables of type $t$ is noted $n^t \in \mathbb{N}$.
The set of type $t$ is noted $\mathcal{X}^t$ and it is such that $\x^t \in \mathcal{X}^t$.

Following the components, a point $\x \in \mathcal{X}$ is defined by a partition by types, and the domain is expressed as Cartesian products of the sets by types, such that
\begin{equation*}
    \x \coloneq \left(
    \x^{\cat},
    \x^{\integer}, \x^{\continuous}\right)
    \in \mathcal{X} \coloneq \mathcal{X}^{\cat} \times \mathcal{X}^{\integer} \times \mathcal{X}^{\continuous}.
\end{equation*}

The categorical variables take \textit{categories} representing qualitative values.
For $i \in I^{\cat}$, the $i$-th categorical variable $x_i^{\cat}$ takes values in the set $\mathcal{X}_i^{\cat} \coloneq \{c_1, c_2, \ldots, c_{\ell_i} \}$, where $c_j$ is a category with $j \in \{1,2,\ldots, \ell_i\}$ and $\ell_i$ is the number of categories of the variable.
For example, the color for an object could be taking a category in the set $\{\text{\sc{red, blue, green}}\}$.
%
{\blr
These variables are inherently difficult to treat because their sets are not ordered and standard distance functions, such as the Euclidean distance, are not suitable for measuring differences between categories~\cite{G-2022-11, hastie01statisticallearning}.
}
The categorical set $\mathcal{X}^{\cat}$ is assumed to be finite.
The total number of categorical components is \smash{$\left| \mathcal{X}^{\cat} \right| = \prod_{i=1}^{n^{\cat}} \ell_i$}.

The integer and continuous variables take quantitative values and they belong to ordered sets in which standard distance functions are well-defined. 
These variables are typically much easier to optimize than the categorical ones.
For convenience, the integer and continuous variables may be regrouped in a quantitative component such that $\x^{\quant} \coloneq \left( \x^{\integer}, \x^{\continuous} \right) \in \mathcal{X}^{\quant} \coloneq \mathcal{X}^{\integer} \times \mathcal{X}^{\continuous}$.

Note that this study does not consider \textit{meta variables}, which influence the inclusion or bounds of other variables~\cite{G-2022-11, G-2024-33}.
The problems have fixed dimensions and bounds.

\subsection{Research gap, objectives and contributions}
\label{sec:objectives_organization}

Bayesian optimization (BO) frameworks have been successfully applied to mixed-variable blackbox problems~\cite{ZhTaChAp2020,PeBrBaTaGu2019, QiWuJe2008, BaDiMoLeSa2023}.
However, BO does not provide local optimization mechanisms, convergence guarantees, or principled stopping criteria.
Several metaheuristics rely on operations based on Gower-type distances that reduce categorical differences to mismatches between categories.
Such generic measures are not problem-informed and can lead to a large number of unnecessary evaluations.

There are other mixed-variable methods, such as \catmads and \mvmads, that handle categorical variables with neighborhoods~\cite{catmads, Abra02a, AACW09a, LuPiSc05a}.
%
In these methods, neighborhoods are used to generate candidate solutions by modifying the variables of an incumbent.
The main advantage of these neighborhoods is that they introduce a local mechanism that leads to theoretical guarantees for discrete variables.
The performance of neighborhood-based approaches depends highly on the quality of the neighborhoods.
In \mvmads, neighborhoods are constructed with user-defined rules, which can be difficult to establish in a blackbox context.
In practice, these neighborhoods are primarily intended to handle discrete variables, but they may also modify the continuous variables of an incumbent.
The theoretical results of \mvmads assume that the points in the user-defined neighborhoods belong to the feasible set $\Omega$.
This constitutes an important methodological restriction for handling constraints and does not capitalize on information from relaxable constraints.
\catmads focuses instead on neighborhoods strictly for categorical variables and constructs them using categorical distances defined solely from the objective function:
The constraint functions are thus not considered.
This can be problematic for problems in which the behavior of the constraint functions changes drastically with respect to the categorical variables.
\Cref{fig:intro_big_picture} illustrates a simple example with one categorical variable
$x^{\cat} \in \{\text{\sc{purple, red, green}}\}$ and two continuous variables.
Each plane corresponds to a category, and the colored regions denote feasibility.
Each plane also shows some level curves of the objective function.

\begin{figure}[htb!]
\centering

\begin{subfigure}[b]{0.45\textwidth}
  \centering
  \scalebox{1.25}{\newcommand{\curveGreen}{
      (-2.2,-1)
      .. (0.2,-1.5)
      .. (0.8,-1.7)
      .. (2,2)
      .. (0,1)
}

\newcommand{\curveBlue}{(-1,0) .. (0.5,3) .. (2,2) .. (0,1)} 

\newcommand{\curveRed}{
  (-2,-0.8)
  .. (0.2,-1.5)
  .. (1.0,-0.8)
  .. (1.5,2)
  .. (0.5,1.2)
  .. (0,1)
}

\begin{tikzpicture}[use Hobby shortcut]
\definecolor{myblue}{RGB}{70,130,255} 
\definecolor{mypurple}{RGB}{186,85,211}
\definecolor{myred}{RGB}{240,60,60}
\definecolor{myorangegreen}{RGB}{210,180,80} 
\definecolor{mygreen}{RGB}{80,180,80}

\definecolor{planeLight}{RGB}{238,238,238}
\definecolor{planeMid}{RGB}{210,210,210}
\definecolor{contourDark}{RGB}{45,45,45}

\tikzset{
  contourLabel/.style={
    font=\scriptsize,
    text=contourDark
  }
}

\def\xA{-0.75} \def\xB{1}
\def\xC{-0.5}  \def\xD{1.25}
\def\yBottom{-0.5}
\def\yTop{1.0}

\def\QuarterPlane{(\xA,\yBottom) -- (\xB,\yBottom) -- (\xD,\yTop) -- (\xC,\yTop) -- cycle;}

\def\Gap{1.6}



\def\LblGreenOneX{-0.55} \def\LblGreenOneY{-0.2} \def\LblGreenOneVal{$0$}
\def\LblGreenTwoX{0.15}  \def\LblGreenTwoY{0.3} \def\LblGreenTwoVal{$3$}
\def\LblGreenTreX{0.75}  \def\LblGreenTreY{0.7} \def\LblGreenTreVal{$6$}

\def\LblTwoTwoX{-0.1} \def\LblTwoTwoY{0.5} \def\LblTwoTwoVal{$4$}
\def\LblTwoOneX{-0.25} \def\LblTwoOneY{0.7} \def\LblTwoOneVal{$6$}
\def\LblTwoTreX{-0.4} \def\LblTwoTreY{0.85} \def\LblTwoTreVal{$8$}

\def\LblThreeTwoX{0.25} \def\LblThreeTwoY{0.5} \def\LblThreeTwoVal{$2$}
\def\LblThreeOneX{0.7} \def\LblThreeOneY{0.75} \def\LblThreeOneVal{$4$}
\def\LblThreeTreX{1} \def\LblThreeTreY{0.85} \def\LblThreeTreVal{$6$}

\def\xStarX{0.62}     
\def\xStarY{0.05}     

\newcommand{\LevelCurvesSimilarBase}{%
  \pgfmathsetmacro{\cx}{0.29/1.75 -0.1}
  \pgfmathsetmacro{\cy}{0.35 -0.1}

  \draw[contourDark, opacity=0.72, line width=0.50pt] (\cx,\cy) ellipse [x radius=0.55, y radius=0.3];
  \draw[contourDark, opacity=0.72, line width=0.50pt] (\cx,\cy) ellipse [x radius=1.10, y radius=0.52];
  \draw[contourDark, opacity=0.72, line width=0.50pt] (\cx,\cy) ellipse [x radius=1.85, y radius=0.7];
}

\newcommand{\LevelCurvesPlaneOne}{%
  \begin{scope}[rotate=5]
    \LevelCurvesSimilarBase
  \end{scope}
}

\newcommand{\LevelCurvesPlaneTwo}{%
  \pgfmathsetmacro{\cx}{0.29/1.75 -0.1}
  \pgfmathsetmacro{\cy}{0.35 -0.1}
  \begin{scope}[shift={(\cx,\cy)}, rotate=20, shift={(-\cx,-\cy)}]
    \LevelCurvesSimilarBase
  \end{scope}
}

\newcommand{\LevelCurvesGreen}{%
  \draw[contourDark, opacity=0.72, line width=0.50pt] (-1.80,-0.2) ellipse [x radius=3.30, y radius=1.40];
  \draw[contourDark, opacity=0.72, line width=0.50pt] (-2.40,-0.85) ellipse [x radius=3.60, y radius=1.55];
  \draw[contourDark, opacity=0.72, line width=0.50pt] (-2.85,-1.45) ellipse [x radius=3.60, y radius=1.65];
}

\foreach \i in {0,...,2} {
  \begin{scope}[shift={(0, \i*\Gap)}, xscale=1.75, yscale=1]
    \ifnum\i=0
      \fill[planeMid, opacity=0.35] \QuarterPlane;
    \else
      \fill[planeLight, opacity=0.70] \QuarterPlane;
    \fi
  \end{scope}
}

\begin{scope}[shift={(1.35, 0.9)}]
  \begin{scope}[shift={(0, 2*\Gap)}, xscale=1.25, yscale=0.7]
    \begin{scope}[scale=0.2, rotate=-30, shift={(-0., -1)}]
      \draw[use Hobby shortcut, closed, fill=mypurple, draw=none, opacity=0.6] \curveRed;
    \end{scope}
  \end{scope}
\end{scope}

\begin{scope}[shift={(0, 1*\Gap)}, xscale=1.75, yscale=1]
  \begin{scope}[scale=0.4, rotate=30, shift={(0.75, 0.25)}]
    \fill[myred, opacity=0.6, draw=none]
      (0, 0) ellipse [x radius=2, y radius=0.8];
  \end{scope}
\end{scope}

\begin{scope}[shift={(0, 0*\Gap)}, xscale=1.75, yscale=1]
  \begin{scope}[scale=0.4, rotate=30, shift={(0.75, 0.25)}]
    \fill[mygreen, opacity=0.6, draw=none]
      (0, 0) ellipse [x radius=1.5, y radius=0.6];
  \end{scope}
\end{scope}

\foreach \i in {0,...,2} {
  \begin{scope}[shift={(0, \i*\Gap)}, xscale=1.75, yscale=1]
    \clip \QuarterPlane;

    \ifnum\i=0
      \LevelCurvesGreen
      \node[contourLabel] at (\LblGreenOneX,\LblGreenOneY) {\LblGreenOneVal};
      \node[contourLabel] at (\LblGreenTwoX,\LblGreenTwoY) {\LblGreenTwoVal};
      \node[contourLabel] at (\LblGreenTreX,\LblGreenTreY) {\LblGreenTreVal};

    \else\ifnum\i=1
      \LevelCurvesPlaneTwo
      \node[contourLabel] at (\LblTwoOneX,\LblTwoOneY) {\LblTwoOneVal};
      \node[contourLabel] at (\LblTwoTwoX,\LblTwoTwoY) {\LblTwoTwoVal};
      \node[contourLabel] at (\LblTwoTreX,\LblTwoTreY) {\LblTwoTreVal};

    \else
      \LevelCurvesPlaneOne
      \node[contourLabel] at (\LblThreeOneX,\LblThreeOneY) {\LblThreeOneVal};
      \node[contourLabel] at (\LblThreeTwoX,\LblThreeTwoY) {\LblThreeTwoVal};
      \node[contourLabel] at (\LblThreeTreX,\LblThreeTreY) {\LblThreeTreVal};
    \fi\fi

  \end{scope}
}

\begin{scope}[shift={(0-0.375, 3.4-0.3-\Gap)}, xscale=1.5, yscale=1]
  \fill[black] (0.29, 0.351) circle (2pt)
    node[above, xshift=15pt, yshift=-5pt] {$\x_{(k)}$};

  \fill[black] (0.29, 1.95) circle (2pt);
  \draw[dashed,-{Latex[length=1.25mm, width=2mm]},thick] (0.29, 0.351) -- (0.29, 1.9);
\end{scope}


\begin{scope}[shift={(0,0*\Gap)}, xscale=1.5, yscale=1]
  \fill[black] ({(\xStarX-1.175)/1.75}, \xStarY-0.15) circle (2pt)
    node[above left] {$\x_\star$};
\end{scope}

\end{tikzpicture}}
  \caption{Categorical neighborhood relying only on objective information.}
  \label{subfig:intro_big_picture_1}
\end{subfigure}\hfill
\begin{subfigure}[b]{0.45\textwidth}
  \centering
  \scalebox{1.25}{\newcommand{\curveGreen}{
      (-2.2,-1)
      .. (0.2,-1.5)
      .. (0.8,-1.7)
      .. (2,2)
      .. (0,1)
}

\newcommand{\curveBlue}{(-1,0) .. (0.5,3) .. (2,2) .. (0,1)} 

\newcommand{\curveRed}{
  (-2,-0.8)
  .. (0.2,-1.5)
  .. (1.0,-0.8)
  .. (1.5,2)
  .. (0.5,1.2)
  .. (0,1)
}

\begin{tikzpicture}[use Hobby shortcut]
\definecolor{myblue}{RGB}{70,130,255} 
\definecolor{mypurple}{RGB}{186,85,211}
\definecolor{myred}{RGB}{240,60,60}
\definecolor{myorangegreen}{RGB}{210,180,80} 
\definecolor{mygreen}{RGB}{80,180,80}

\definecolor{planeLight}{RGB}{238,238,238}
\definecolor{planeMid}{RGB}{210,210,210}
\definecolor{contourDark}{RGB}{45,45,45}

\tikzset{
  contourLabel/.style={
    font=\scriptsize,
    text=contourDark
  }
}

\def\xA{-0.75} \def\xB{1}
\def\xC{-0.5}  \def\xD{1.25}
\def\yBottom{-0.5}
\def\yTop{1.0}

\def\QuarterPlane{(\xA,\yBottom) -- (\xB,\yBottom) -- (\xD,\yTop) -- (\xC,\yTop) -- cycle;}

\def\Gap{1.6}



\def\LblGreenOneX{-0.55} \def\LblGreenOneY{-0.2} \def\LblGreenOneVal{$0$}
\def\LblGreenTwoX{0.15}  \def\LblGreenTwoY{0.3} \def\LblGreenTwoVal{$3$}
\def\LblGreenTreX{0.75}  \def\LblGreenTreY{0.7} \def\LblGreenTreVal{$6$}

\def\LblTwoTwoX{-0.1} \def\LblTwoTwoY{0.5} \def\LblTwoTwoVal{$4$}
\def\LblTwoOneX{-0.25} \def\LblTwoOneY{0.7} \def\LblTwoOneVal{$6$}
\def\LblTwoTreX{-0.4} \def\LblTwoTreY{0.85} \def\LblTwoTreVal{$8$}

\def\LblThreeTwoX{0.25} \def\LblThreeTwoY{0.5} \def\LblThreeTwoVal{$2$}
\def\LblThreeOneX{0.7} \def\LblThreeOneY{0.75} \def\LblThreeOneVal{$4$}
\def\LblThreeTreX{1} \def\LblThreeTreY{0.85} \def\LblThreeTreVal{$6$}

\def\xStarX{0.62}     
\def\xStarY{0.05}     

\newcommand{\LevelCurvesSimilarBase}{%
  \pgfmathsetmacro{\cx}{0.29/1.75 -0.1}
  \pgfmathsetmacro{\cy}{0.35 -0.1}

  \draw[contourDark, opacity=0.72, line width=0.50pt] (\cx,\cy) ellipse [x radius=0.55, y radius=0.3];
  \draw[contourDark, opacity=0.72, line width=0.50pt] (\cx,\cy) ellipse [x radius=1.10, y radius=0.52];
  \draw[contourDark, opacity=0.72, line width=0.50pt] (\cx,\cy) ellipse [x radius=1.85, y radius=0.7];
}

\newcommand{\LevelCurvesPlaneOne}{%
  \begin{scope}[rotate=5]
    \LevelCurvesSimilarBase
  \end{scope}
}

\newcommand{\LevelCurvesPlaneTwo}{%
  \pgfmathsetmacro{\cx}{0.29/1.75 -0.1}
  \pgfmathsetmacro{\cy}{0.35 -0.1}
  \begin{scope}[shift={(\cx,\cy)}, rotate=20, shift={(-\cx,-\cy)}]
    \LevelCurvesSimilarBase
  \end{scope}
}

\newcommand{\LevelCurvesGreen}{%
  \draw[contourDark, opacity=0.72, line width=0.50pt] (-1.80,-0.2) ellipse [x radius=3.30, y radius=1.40];
  \draw[contourDark, opacity=0.72, line width=0.50pt] (-2.40,-0.85) ellipse [x radius=3.60, y radius=1.55];
  \draw[contourDark, opacity=0.72, line width=0.50pt] (-2.85,-1.45) ellipse [x radius=3.60, y radius=1.65];
}

\foreach \i in {0,...,2} {
  \begin{scope}[shift={(0, \i*\Gap)}, xscale=1.75, yscale=1]
    \ifnum\i=0
      \fill[planeMid, opacity=0.35] \QuarterPlane;
    \else
      \fill[planeLight, opacity=0.70] \QuarterPlane;
    \fi
  \end{scope}
}

\begin{scope}[shift={(1.35, 0.9)}]
  \begin{scope}[shift={(0, 2*\Gap)}, xscale=1.5, yscale=0.75]
    \begin{scope}[scale=0.2, rotate=-30, shift={(-0., -1)}]
      \draw[use Hobby shortcut, closed, fill=mypurple, draw=none, opacity=0.6] \curveRed;
    \end{scope}
  \end{scope}
\end{scope}

\begin{scope}[shift={(0, 1*\Gap)}, xscale=1.75, yscale=1]
  \begin{scope}[scale=0.4, rotate=30, shift={(0.75, 0.25)}]
    \fill[myred, opacity=0.6, draw=none]
      (0, 0) ellipse [x radius=2, y radius=0.8];
  \end{scope}
\end{scope}

\begin{scope}[shift={(0, 0*\Gap)}, xscale=1.75, yscale=1]
  \begin{scope}[scale=0.4, rotate=30, shift={(0.75, 0.25)}]
    \fill[mygreen, opacity=0.6, draw=none]
      (0, 0) ellipse [x radius=1.5, y radius=0.6];
  \end{scope}
\end{scope}

\foreach \i in {0,...,2} {
  \begin{scope}[shift={(0, \i*\Gap)}, xscale=1.75, yscale=1]
    \clip \QuarterPlane;

    \ifnum\i=0
      \LevelCurvesGreen
      \node[contourLabel] at (\LblGreenOneX,\LblGreenOneY) {\LblGreenOneVal};
      \node[contourLabel] at (\LblGreenTwoX,\LblGreenTwoY) {\LblGreenTwoVal};
      \node[contourLabel] at (\LblGreenTreX,\LblGreenTreY) {\LblGreenTreVal};

    \else\ifnum\i=1
      \LevelCurvesPlaneTwo
      \node[contourLabel] at (\LblTwoOneX,\LblTwoOneY) {\LblTwoOneVal};
      \node[contourLabel] at (\LblTwoTwoX,\LblTwoTwoY) {\LblTwoTwoVal};
      \node[contourLabel] at (\LblTwoTreX,\LblTwoTreY) {\LblTwoTreVal};

    \else
      \LevelCurvesPlaneOne
      \node[contourLabel] at (\LblThreeOneX,\LblThreeOneY) {\LblThreeOneVal};
      \node[contourLabel] at (\LblThreeTwoX,\LblThreeTwoY) {\LblThreeTwoVal};
      \node[contourLabel] at (\LblThreeTreX,\LblThreeTreY) {\LblThreeTreVal};
    \fi\fi

  \end{scope}
}

\begin{scope}[shift={(0-0.375, 3.4-0.3-\Gap)}, xscale=1.5, yscale=1]
  \fill[black] (0.29, 0.351) circle (2pt)
    node[above, xshift=15pt, yshift=-5pt] {$\x_{(k)}$};

  \fill[black] (0.29, 1.95-2*\Gap) circle (2pt);
  \draw[dashed,-{Latex[length=1.25mm, width=2mm]},thick] (0.29, 0.351) -- (0.29, 2-2*\Gap);
\end{scope}

\begin{scope}[shift={(0,0*\Gap)}, xscale=1.5, yscale=1]
  \fill[black] ({(\xStarX-1.175)/1.75}, \xStarY-0.15) circle (2pt)
    node[above left] {$\x_\star$};
\end{scope}

\end{tikzpicture}} 
  \caption{Categorical neighborhood integrating objective and constraints information.}
  \label{subfig:intro_big_picture_2}
\end{subfigure}
\caption{A motivating example where the incumbent $\x_{(k)}$ has objective value $f(\x_{(k)})=3$, while the global minimum $\x_{\star}$ has $f(\x_{\star})=1$.
The neighborhood on the left selects {\sc purple} based strictly on proximity with the objective.
The neighborhood on the right selects {\sc green} based on proximity of both the constraints and the objective. 
}
\label{fig:intro_big_picture}
\end{figure}

In \cref{fig:intro_big_picture}, the current solution $\x_{(k)}$ lies in the continuous subspace associated with the categorical component $\x^{\cat}=\text{\sc{red}}$.
It corresponds to a local minimum with respect to the continuous variables.
However, the global minimum $\x_{\star}$ resides in the continuous subspace associated with $x^{\cat}=\text{\sc{green}}$. 
In the setting, reaching this global minimum from $\x_{(k)}$ can be difficult when using a neighborhood-based approach.
From $\x_{(k)}$ in the subspace of $\text{\textsc{red}}$, objective-based or user-defined neighborhoods would likely prioritize the category $\text{\sc{purple}}$.
The restriction of the objective function in the $\text{\sc{red}}$ and $\text{\sc{purple}}$ subspaces appears more similar.
%
Consequently, an objective-based neighborhood may focus on the $\text{\sc{purple}}$ category, although: 1) its associated continuous subspace is highly infeasible and yields poorer feasible values, and 2) the feasible regions in the continuous subspaces of $\text{\sc{red}}$ and $\text{\sc{green}}$ are much more similar.
%
%
This can degrade the performance of neighborhood-based methods and the quality of the local minima with respect to the categorical variables.
To mitigate this issue, one could rely on global mechanisms such as BO or design better-informed neighborhoods.
In the presence of noise or nonsmoothness, BO may fail to properly tackle the continuous variables and by consequence the mixed-variable problem itself.
More local and robust mechanisms may therefore be needed to ensure incremental progress.

This motivates the development of local neighborhoods that are better-informed and robust.
The research gap addressed is the automatic construction of problem-specific categorical neighborhoods that incorporate information from both the objective and the constraint functions in a blackbox setting.
Note that, hereafter, categorical neighborhoods are referred to simply as neighborhoods.
The work does not aim to introduce a new optimization method, but rather to provide additional structural information that facilitates handling categorical variables.
These better-informed neighborhoods provide a systematic structure for categorical variables, improve interpretability, and are particularly well suited for blackbox constrained optimization.
They enhance the quality of local minima with respect to categorical variables and improve the efficiency {\blr of} neighborhood-based methods.
The number of evaluations for optimizing categorical variables may decrease, since each categorical neighbor is selected in a more informed manner.

The present work is structured through four objectives, representing the main contributions of the work. 
The first objective is to construct a problem-specific categorical pseudo-distance using kernel-based surrogate models, such as Gaussian Processes (GPs), so that categorical components with similar predicted objective values are considered closer.
%
%
The second objective is to define a categorical pseudo-distance that captures similarity with respect to the constraint functions by aggregating surrogate models of the constraints.
Building on these two pseudo-distances, the third objective is to construct categorical neighborhoods balancing similarity between the objective and the constraints.
The first three objectives are developed independently of any specific optimization framework and for purely categorical domains.
The fourth objective is to generalize the proposed neighborhoods to mixed-variable domains.
As a proof of concept for the fourth objective, the neighborhoods are used in \catmads.

\subsection{Related work}
\label{sec:intro_related}

%
In~\cite{catmads}, a categorical distance is constructed using a basic interpolation of the objective function.
Neighborhoods are characterized with this single distance function and they are constructed without accounting {\blr for} the constraint functions, which may lead to neighborhoods containing elements that are highly dissimilar in terms of feasibility.
In~\cite{AACW09a,LuPiSc05a}, categorical neighborhoods are specified manually by the user, which requires prior knowledge about the problem structure, which may be ambiguous in a blackbox context.
The primary focus of \catmads and \mvmads lies in the optimization process, while neighborhood design remains secondary.
In contrast, the present work emphasizes the systematic construction of categorical neighborhoods, independently of any particular optimization framework.

Several metaheuristic methods have been proposed to address categorical variables in mixed-variable optimization.
These methods construct candidate points with operations transforming categorical variables of existing solutions.
%
Such operations have similar roles to categorical neighborhoods.
%
CatCMA extends the covariance matrix adaptation evolution strategy (CMA-ES) to problems involving both categorical and continuous variables~\cite{HaSaNoUcSh2024}.
In this approach, categorical variables are modeled by probability vectors over their respective categories, while continuous variables are modeled by a multivariate normal distribution.
A joint probabilistic distribution over categorical and continuous variables is introduced to guide the optimization in mixed-variable domains.
{\blr 
The Non-dominated Sorting Genetic Algorithm-II (NSGA-II)~\cite{DePrAgMe2002}, although originally developed for multi-objective optimization, was extended to mixed-variable problems involving categorical variables.
An implementation proposed in~\cite{BlankDeb2020} uses crossover and mutation operators tailored to each variable type to generate candidate points from a population of solutions.
%
} 
%
%
For more details, the survey~\cite{Ta2024} provides an exhaustive survey for mixed-variable metaheuristics. 

In Bayesian optimization, categorical variables are structured with similarity measures, called \textit{kernels}.
In optimization, kernels are problem-specific similarity measures characterizing surrogate models that guide the optimization.
In the literature, various techniques exist for constructing categorical kernels, most notably \textit{one-hot encoding} and matrix-based approaches.
One-hot encoding assigns a unique binary variable to each category, and they are assigned their corresponding hyperparameter weighting the importance of each category in the kernel~\cite{GaHe2020}.
Matrix-based approaches have shown great results in BO~\cite{PeBrBaTaGu2019, QiWuJe2008, RoPaDeClPeGiWy2020, BaDiMoLeSa2023}. 
%
The elements of these matrices are hyperparameters that model correlations and anti-correlations between different categories~\cite{QiWuJe2008}.
The similarities between categories are directly encoded in matrices.
Other techniques for encoding categorical variables include continuous latent spaces.
In~\cite{ZhTaChAp2020}, each categorical variable is assigned a 2D continuous space, in which its categories are mapped in a way that correlated categories are close, and uncorrelated ones are far. 
The coordinates of the categories can be viewed as hyperparameters representing correlations.  
%

The rest of this document is organized as follows. 
An illustrative optimization problem is introduced in \cref{sec:example} to guide the presentation.
The essential concepts of kernels and GPs are presented in \cref{sec:kernel}.
In \cref{sec:neighborhoods}, the categorical pseudo-distances and neighborhoods are developed.
The novel neighborhoods are adapted and used within the \catmads framework in \cref{sec:direct_search} to provide  proof of concept.
Lastly, some numerical experiments benchmarking \catmads using the new neighborhoods with other solvers are presented in \cref{sec:computational_experiments}.

\section{An illustrative mechanical-part design problem}
\label{sec:example}

To outline the proposed contributions, a mechanical-part design problem is used as an illustrative example in the following sections.
The goal is to minimize the physical strain of a structural piece subject to a budget constraint and an ecological score constraint.
A point represents a design, and it is composed of three categorical variables and two continuous variables.
The categorical variables correspond to the choice of supplier $u \in \{\A, \B\}$, the material $a \in \{\Aluminum, \Steel, \Composite\footnote{\Composite for ``composite''}, \Wood\}$, and the shape $s \in \{\Square, \Circle, \Ellipse\}$.
%
%
The continuous variables are a length $l \in [5,10]$ and a ratio of recycle{\bl d} material $r \in [0,1]$.
The domain and a point are defined as 
\begin{equation*}
  \x = (u, a, s, l, r) \in \mathcal{X}
  = \mathcal{X}^{\cat} \times \mathcal{X}^{\continuous} 
\end{equation*}
where $\mathcal{X}^{\cat} = \{\A, \B\}
     \times \{\Aluminum, \Steel, \Composite, \Wood\}  \times \{\Square, \Circle, \Ellipse\}$, and
$\mathcal{X}^{\continuous} =
      [5,10]
     \times [0,1]$.
The categorical set $\mathcal{X}^{\cat}$ contains $2 \times 4 \times 3=24$ categorical components, and the feasible region $\Omega$ is unknown.
%

%
%

{\blr
The problem is intentionally constructed as a blackbox in accordance with the scope of this work.
The functions are evaluated through iterative numerical routines to a prescribed accuracy and are therefore not analytically available.
%
%
The blackbox is intentionally inexpensive to allow large statistical samples in order to analyze the effect of categorical variables on constraint feasibility.
Moreover, the problem was designed so that the categorical variables strongly affect constraint feasibility.
Note that the optimization experiments in \cref{sec:computational_experiments} nevertheless use small evaluation budgets to remain representative of the expensive blackbox setting considered in this work.
}

%
%
%
%
%
Each categorical component $\x^{\cat} = (u,a,s) \in \mathcal{X}^{\cat}$ induces a different feasible region over the continuous subspace $\mathcal{X}^{\continuous} = [5,10] \times [0,1]$.
Some categories may appear intuitively similar, \textit{e.g.}, the circle and ellipse categories for the cross-section shape.
However, the objective function, and, in particular, the constraint functions, behave very differently across categories.
%
%
%
As a result, constructing categorical neighborhoods without incorporating information from both the objective and the constraints can lead to misleading similarity measures and poor optimization performance.
The difficulty of the problem is outlined in \cref{fig:histograms} that describes three Latin Hypercube Sampling (LHS) experiments and a graph.
\begin{figure}[htb!]
\centering
\begin{subfigure}[t]{0.4\textwidth}
\centering
\vbox to 7cm{
\centering
\scalebox{0.9}{
\begin{tikzpicture}
\begin{axis}[
    ybar,
    bar width=5pt,
    width=\linewidth*1.1,
    height=6cm,
    ymin=0,
    ymax=10.0,
    x=1.05cm,
    enlarge x limits=0.1,
    ytick={0,2.5,5,7.5,10,12.5},
    ylabel={\% feasible},
    ylabel style={xshift=20pt},
    ylabel near ticks,
    symbolic x coords={
        {(B, Wood, Circle)},
        {(A, Wood, Ellipse)},
        {(A, Composite, Square)},
        {(A, Aluminum, Circle)},
        {(B, Steel, Square)},
        {(B, Aluminum, Square)}
    },
    xtick={
        {(B, Wood, Circle)},
        {(A, Wood, Ellipse)},
        {(A, Composite, Square)},
        {(A, Aluminum, Circle)},
        {(B, Steel, Square)},
        {(B, Aluminum, Square)}
    },
    x tick label style={rotate=45, anchor=east, font=\small},
    xticklabels={
        {$(\B, \Wood, \Circle)$},
        {$(\A, \Wood, \Ellipse)$},
        {$(\A, \Composite, \Square)$},
        {$(\A, \Aluminum, \Circle)$},
        {$(\B, \Steel, \Square)$},
        {$(\B, \Aluminum, \Square)$}
    },
    legend style={
        at={(0.5,1.03)},
        anchor=south,
        legend columns=3,
        draw=none,
        column sep=14pt
    },
    legend image code/.code={
        \draw[#1] (0cm,-0.08cm) rectangle (0.18cm,0.18cm);
    },
    nodes near coords=false,
    unbounded coords=jump
]

\addplot+[
    fill=blue!35,
    every node near coord/.append style={draw=none, fill=none, text opacity=0},
    nodes near coords={},
] coordinates {
    ({(B, Wood, Circle)}, 9.2)
    ({(A, Wood, Ellipse)}, 6.4)
    ({(A, Composite, Square)}, 0.4)
    ({(A, Aluminum, Circle)}, nan)
    ({(B, Steel, Square)}, nan)
    ({(B, Aluminum, Square)}, nan)
};

\addplot+[
    fill=red!35,
    every node near coord/.append style={draw=none, fill=none, text opacity=0},
    nodes near coords={},
] coordinates {
    ({(B, Wood, Circle)}, 7.98)
    ({(A, Wood, Ellipse)}, 8.24)
    ({(A, Composite, Square)}, 0.2)
    ({(A, Aluminum, Circle)}, 0.34)
    ({(B, Steel, Square)}, 0.1)
    ({(B, Aluminum, Square)}, nan)
};

\addplot+[
    fill=green!35,
    every node near coord/.append style={draw=none, fill=none, text opacity=0},
    nodes near coords={},
] coordinates {
    ({(B, Wood, Circle)}, 8.09)
    ({(A, Wood, Ellipse)}, 7.82)
    ({(A, Composite, Square)}, 0.22)
    ({(A, Aluminum, Circle)}, 0.25)
    ({(B, Steel, Square)}, 0.26)
    ({(B, Aluminum, Square)}, 0.095)
};

\legend{250, 5000, 100000}

\end{axis}
\end{tikzpicture}
}
\vfill
}
\caption{Feasibility rate per categorical component.}
\label{subfig:histograms}
\end{subfigure}
\hfill
\begin{subfigure}[t]{0.4\textwidth}
\centering
\vbox to 6.5cm{
\centering
\scalebox{0.9}{
\begin{tikzpicture}
\begin{axis}[
    width=1.2\linewidth,
    height=6.2cm,
    xmode=log,
    xmin=0.008, xmax=12,
    xlabel={\% feasible},
    ymin=0, ymax=32,
    ytick={0,6,12,18,24,30},
    ylabel={Best objective value},
    grid=both,
    clip=false,
    minor grid style={gray!20},
    major grid style={gray!40},
]

\addplot+[only marks, mark=*] coordinates {
    (8.09, 28.28)
    (7.82, 24.09)
    (0.27, 9.98)   
    (0.24, 18.26)  
    (0.21, 8.79)   
    (0.01, 13.95)
};

\node[anchor=east, font=\small] at (axis cs:8.09,29.2) {(\B, \Wood, \Circle)};
\node[anchor=east, yshift=8pt, font=\small] at (axis cs:7.82,23.1) {(\A, \Wood, \Ellipse)};
\node[anchor=west, font=\small] at (axis cs:0.27,10.9) {(\B, \Steel, \Square)};
\node[anchor=west, font=\small] at (axis cs:0.24,19.3) {(\A, \Aluminum, \Circle)};
\node[anchor=west, font=\small] at (axis cs:0.21,7.7) {(\A, \Composite, \Square)};
\node[anchor=west, font=\small] at (axis cs:0.01,14.9) {(\B, \Aluminum, \Square)};

\end{axis}
\end{tikzpicture}
}
\vfill
}
\caption{Best objective function value vs.~\% feasible with $100{,}000$ points per categorical component.}
\label{subfig:scatter3}
\end{subfigure}

\caption{Statistical analysis of feasibility per categorical component in
{\blr the mechanical-part design problem.}}
\label{fig:histograms}
\end{figure}

%
In the histogram of \cref{subfig:histograms}, LHS is performed independently within each categorical component.
Each bar represents the percentage of feasible points obtained for a given number of samples.
%
%
%
Only the categorical {\blr components} for which at least one feasible point was sampled {\blr are} presented.
%
%
As the number of points per categorical component increases, more components with feasible points are discovered.
%
The categorical component $\left(\A, \Aluminum, \Circle\right)$ does not have any feasible point with $250$ samples (there is no blue bin), but it does with $5{,}000$ and $100{,}000$ samples. 
The histogram suggests that only six categorical components may be capable of achieving feasibility.

The graph in \cref{subfig:scatter3} shows the best feasible objective value obtained by performing an LHS with $100{,}000$ points for each categorical component.
Overall, best objective values are associated with smaller percentages of feasibility.

\section{Mixed-variable kernels}
\label{sec:kernel}

{\blr
This study uses similarity measures derived from available data to induce categorical pseudo-distances, including those derived from kernels.
}
A similarity measure quantifies how close or related two data points are in the decision space. 
A kernel is a particular type of similarity measure.
The data points of Problem~\eqref{eq:formulation} lie in the domain~$\mathcal{X}$, and the kernel is defined on this domain.
%
%
{\blr A categorical kernel is then extracted from the kernel defined on $\mathcal{X}$ and used to construct categorical neighborhoods.}
In practice, kernels can compare mixed-variable points, enabling the construction of interpolation models that rely exclusively on similarity measures, such as mixed-variable Gaussian processes.

Formally, a kernel $\kappa: \mathcal{X} \times \mathcal{X} \to \mathbb{R}$ is a symmetric similarity measure, such that 
\begin{enumerate}
    \item $\kappa(\x,\y) = \kappa(\y,\x)$ for any $\x,\y \in \mathcal{X}$, and 

    \item for any finite set of points $\{\x_{(1)}, \x_{(2)}, \ldots, \x_{(p)}\}$ with $p \in \mathbb{N}$,  the symmetric $p \times p$ matrix $[K]_{i,j}= \kappa \left( \x_{(i)}, \x_{(j)}\right)$ for $i,j \in \{1,2,\ldots, p\}$ is positive semi-definite.
\end{enumerate}
Such {\blr a} kernel $\kappa$ is said to be positive semi-definite.
%
%
The kernel is constructed following the approach of the classic textbook~\cite{RaWi06}. 
For $t \in \{\cat, \quant\}$ and $i \in I^t$,
each variable $x_i^t \in \mathcal{X}_i^t$ is assigned a one-dimensional kernel $\kappa_{i}^t:\mathcal{X}_i^t \times \mathcal{X}_i^t \to \mathbb{R}$ of the appropriate variable type.
%
The one-dimensional kernels are combined with multiplications that conserve symmetry
{\blr and positive semi-definiteness.}

Common kernels for quantitative variables include polynomial, Matern or Gaussian kernels, \textit{e.g.}, see~\cite{RaWi06}.
Typically, the quantitative kernel can be built directly as {\blr a} product of one-dimensional Gaussian kernels, \textit{i.e.},
%
    %
%
\begin{equation*}
    \kappa^{\quant} \left( \x^{\quant}, \y^{\quant} \right) \coloneq 
    \prod_{i=1}^{n^{\quant}} \exp \left( -\theta_i^{\quant} \left(x_i^{\quant} - y_i^{\quant}\right)^2 \right),
\end{equation*}
where $\bm{\theta}^{\quant} \coloneq \left(\theta_1^{\quant}, \theta_2^{\quant}, \ldots, \theta_{n^{\quant}}^{\quant} \right) \in \mathbb{R}_+^{n^{\quant}}$ is the vector of hyperparameters of the quantitative kernel.
The quantitative kernel and the following ones are all parameterized by hyperparameters.
For readability, hyperparameters are not explicitly noted in the arguments of kernels.

\subsection{Categorical kernels}
\label{sec:categorical_kernel}

As mentioned in \cref{sec:intro_related}, there are currently two main state-of-the-art approaches for constructing categorical kernels: one-hot encoding and matrix-based methods.
The matrix-based approach encodes similarity measures between categories directly into matrices.
{\blr
More precisely, each categorical variable $x_i^{\cat} \in \mathcal{X}_i^{\cat}$ is assigned a symmetric positive semi-definite matrix
$T_i^{\cat} \in \mathbb{R}^{\ell_i^{\cat} \times \ell_i^{\cat}}$,
where $\ell_i^{\cat}$ is the number of categories of $x_i^{\cat}$.
The matrix is constructed as
$T_i^{\cat} \coloneq L_i L_i^{\top}$,
where $L_i$ is a lower-triangular matrix parameterized by a hypersphere decomposition (lengths and angles)~\cite{BrMeRa2007, JaRe1999}.
The elements of $L_i$ contain hyperparameters of the categorical kernel, whereas the elements of $T_i^{\cat}$ represent similarities between categories of the corresponding categorical variable.
For example, in the mechanical-part design problem, the similarity matrix associated with the material choice
$a \in \{\Aluminum \text{ (A)}, \Steel \text{ (S)}, \Composite \text{ (C)}, \Wood \text{ (W)}\}$ can be expressed as
\begin{equation*}
T^{\cat}_{a}=
    \begin{bmatrix}
    \vartheta_{\text{A},\text{A}} &   &   & \\
    \vartheta_{\text{S},\text{A}} & \vartheta_{\text{S},\text{S}} & & \\
    \vartheta_{\text{C},\text{A}} & \vartheta_{\text{C},\text{S}} & \vartheta_{\text{C},\text{C}} & \\
    \vartheta_{\text{W},\text{A}} & \vartheta_{\text{W},\text{S}} & \vartheta_{\text{W},\text{C}} & \vartheta_{\text{W},\text{W}}
    \end{bmatrix}.
\end{equation*}
The values of $\vartheta_{i,j}$ are implicitly determined through the hyperparameters of the lower-triangular matrix $L_a$, such that
$T_a^{\cat}=L_aL_a^{\top}$.
The diagonal terms $\vartheta_{i,i}$ represent the similarity of each category with itself.
Since $T_a^{\cat}$ is positive semi-definite, these diagonal terms satisfy
$\vartheta_{i,i}\geq 0$.
In the homoscedastic parametrization, these self-similarities are fixed to
$\vartheta_{i,i}=1$ and can therefore be interpreted as auto-correlations.
In this case, the off-diagonal terms satisfy
$\vartheta_{i,j}\in[-1,1]$.
In the heteroscedastic parametrization, the diagonal terms are allowed to vary while remaining nonnegative.
More generally, positive semi-definiteness implies
$|\vartheta_{i,j}| \leq \sqrt{\vartheta_{i,i}\vartheta_{j,j}}$.
}

For the matrix-based approach, the categorical kernel is derived from the categorical matrices $T_1^{\cat}, T_2^{\cat}, \ldots, T_{n^{\cat}}^{\cat}$.
Essentially, two points $\x$ and $\y$ are compared variable-wise, with each matrix $T_i^{\cat}$ 
providing the similarity measure between categories $x_i^{\cat}$ and $y_i^{\cat}$.
These similarity measures are then multiplied.
%
The number of hyperparameters associated with the categorical variable $x_i^{\cat}$ depends on the number of categories $\ell_i$ and the chosen parametrization of $L_i$.
While more hyperparameters increase modeling flexibility, they also imply greater computational costs:
see~\cite{BaDiMoLeSa2023} for a detailed presentation of parametrizations and the hypersphere decomposition.
The matrix-based approach typically involves many hyperparameters.
Alternative approaches with fewer hyperparameters are available.
The common \textit{one-hot encoding} technique assigns a unique binary variable to each category. 
The categorical component $\x^{\cat} \in \mathcal{X}^{\cat}$ is represented by a total of $ \sum_{i=1}^{n^{\cat}} \ell_i$
binary variables, where $\ell_i \in \mathbb{N}$ is the number of categories of the $i$-th categorical variable $x_i^{\cat}$.
For each categorical variable, exactly one of its binary variables is assigned the value $1$, and the others are set to $0$.
{\blr In the mechanical-part design problem,}
recall that the categorical variables are $u \in \{\A, \B\}$, $a \in \{\Aluminum, \Steel, \Composite, \Wood \}$ and
$ s \in \{ \Square, \Circle, \Ellipse \}$.
The categorical component $(\B, \Wood, \Circle)$ would be represented by $\left( \left(0,1 \right), \left(0,0,0,1\right), \left(0,1,0 \right) \right)$, whereas $(\A, \Steel, \Ellipse)$ would be represented by $\left( \left(1,0 \right), \left(0,1,0,0\right), \left(0,0,1 \right) \right)$.
With this approach, each binary variable can be assigned a hyperparameter.
%
Then, a categorical kernel can be defined via Gaussian kernels as follows
\begin{equation}
    \kappa^{\cat} \left( \x^{\cat}, \y^{\cat} \right) \coloneq 
    \prod_{i=1}^{n^{\cat}} \exp \left( - \left\Vert \bm{E}_i\left( x_i^{\cat} \right) - \bm{E}_i\left( y_i^{\cat} \right) \right\Vert^2 _{\bm{\theta}_i^{\cat}} \right)
\label{eq:categorical_kernel}
\end{equation}
where $\bm{E}_i\left( x_i^{\cat} \right) \in \left\{ \z \in \{0,1\}^{\ell_i} \, : \, \mathbf{1}^{\top} \z = 1 \right\}$
is the \textit{one-hot} representation of the $i$-th categorical variable $x_i^{\cat}$, and $\bm{\theta}_i^{\cat} \in \mathbb{R}_+^{\ell_i}$ is the vector of hyperparameters of $x_i^{\cat}$,
and {\blr $\left\Vert \bm{u} - \bm{v} \right\Vert_{\bm{w}}^2 \coloneq (\bm{u}-\bm{v})^{\top} \text{diag} \left(\bm{w}\right) (\bm{u}-\bm{v})$ is the squared weighted Euclidean norm.
%
%
%
In~\eqref{eq:categorical_kernel}, the squared norm represents squared distances between binary vectors, where the weights are the hyperparameters of the binary variables.}

%
%
%

\subsection{Adjusting the kernels with mixed-variable GPs}
\label{sec:HPs_and_GPs}

%
{\blr Once appropriate quantitative and categorical kernels are defined, the kernel $\kappa:\mathcal{X} \times \mathcal{X} \to \mathbb{R}$ can be constructed as the product of the categorical and quantitative kernels:
%
\begin{equation*}
    \kappa \left(\x,\y \right) \coloneq \kappa^{\cat} \left(\x^{\cat}, \y^{\cat} \right) \,
    \kappa^{\quant} \left(\x^{\quant}, \y^{\quant} \right),
\end{equation*}
which is parameterized by a vector of hyperparameters $\bm{\theta}\coloneq \left(\bm{\theta}^{\cat}, \bm{\theta}^{\quant} \right)$.
}

To adjust the hyperparameters of the kernel, a regression model linking input data points to their corresponding images is required.
In BO, GPs are commonly used as surrogate models for objective and constraint functions, as they can represent a wide variety of functions, including mixed-variable ones.
%

%
For readability, the construction of GP surrogates is presented for the objective function.
%
%
This surrogate is noted $\tilde{f}$ and it is constructed with a set of points $\mathbb{X}\coloneq \{\x_{(1)}, \x_{(2)},\ldots, \x_{(p)}\}$ and a vector of corresponding images $\bm{f} =(f(\x_{(1)}), f(\x_{(2)}), \ldots, f(\x_{(p)}))$.
{\blr
Only evaluations with finite output values are used to fit the Gaussian process.
Infinite values cannot be used in the interpolation model and are therefore excluded.
}
This procedure yields a prediction function $\hat{f}:\mathcal{X}\to\mathbb{R}$
and a variance function 
$\hat{\sigma}^2:\mathcal{X} \to \mathbb{R}_+$ 
that, in the noiseless setting, can be expressed as in~\cite{RaWi06}:
\begin{align*}
\begin{split}
    \hat{f}(\x) &\coloneq \bm{\kappa}(\x)^{\top} K^{-1}  \bm{f} , \\
    {\hat{\sigma}}^2(\x) &\coloneq  \kappa(\x,\x) - \bm{\kappa}(\x)^{\top} K^{-1} \bm{\kappa}(\x) \geq 0,
\end{split}
\end{align*}
where $\kappa:\mathcal{X} \times \mathcal{X} \to \mathbb{R}$ is the kernel, $\bm{\kappa}(\x)\coloneq\left( \kappa(\x, \x_{(1)}), \kappa(\x, \x_{(2)}), \ldots, \kappa(\x, \x_{(p)})  \right) \in \mathbb{R}^{p}$
is the kernel vector comparing an input point $\x \in \mathcal{X}$ and the points in $\mathbb{X}$, and $K \in \mathbb{R}^{p \times p}$ is the kernel matrix that compares all points in $\mathbb{X}$ such that 
$[K]_{i,j}\coloneq \kappa \left(\x_{(i)}, \x_{(j)} \right)$ for $i,j \in \{1,2,\ldots, p\}$.

The kernel characterizes the GP by quantifying how correlated or similar input points are, \textit{i.e.}, how close their output should be.
Concretely, it controls the general behavior of the fit, \textit{e.g.}, the smoothness in a continuous setting.
In the mixed-variable setting, defining an appropriate mixed-variable kernel is key to constructing a GP surrogate that can model the functions of interest.
In fact, the GP that best fits the available data is found by adjusting the hyperparameters of the kernel.
They are typically obtained by maximizing the marginal likelihood of the GP~\cite{RaWi06}:
\begin{equation}
\bm{\theta}_{\star} \in \argmax_{\bm{\theta}} \log  \mathbb{P} \left[ \bm{f} \mid \mathbb{X}, \bm{\theta} \right]
\quad \mbox{ where }\quad 
    \log  \mathbb{P} \left[ \bm{f} \mid \mathbb{X}, \bm{\theta} \right]  = - \frac{1}{2} \bm{f}^{\top} K^{-1} \bm{f} - \frac{1}{2} \log  \det \left( K\right)  - \frac{p}{2} \log \left( 2\pi \right).
\label{eq:log_likelihood}
\end{equation}
{\blr 
The marginal likelihood is continuously differentiable.
It can be optimized with standard continuous solvers using a multistart strategy.
Several starting points are sampled in the hyperparameter space and used to initialize the solver.
The construction of the GP also involves operations with matrix $K$, which may become computationally expensive or numerically unstable when $K$ is large or ill-conditioned.
Typically, the data is normalized and a small nugget term is added to the diagonal terms of $K$ to improve numerical stability.
A limitation of using GPs is that these numerical issues may become more significant as the number of sampled points increases.
If $K$ remains ill-conditioned, points that are too close may need to be removed or the nugget term may need to be increased.
}

{\blr After this step, the GP is properly constructed and the categorical kernel is finely adjusted to the problem data.
The categorical kernel is used in the next section.}

The GP surrogates of the constraint functions are constructed using the same procedure as for the objective function.
Specifically, for $j \in J$, each constraint function $g_j$ is provided a probabilistic surrogate $\tilde{g}_j$ with a predictive mean function $\hat{g}_j:\mathcal{X} \to \mathbb{R}$ and variance function $\hat{\sigma}_j^2:\mathcal{X} \to \mathbb{R}_+$, as well as its own kernel and hyperparameters.

\section{Surrogate-based categorical neighborhoods}
\label{sec:neighborhoods}

%

{\blr Kernels and surrogate models are used to construct neighborhoods that induce a structure over the categorical variables.}
The proposed neighborhoods are defined using two categorical pseudo-distances, one associated with the objective function and another one associated with the constraint functions.
Proximity between categorical variables simultaneously considers similarities in the objective and constraint functions.

This section focuses exclusively on categorical variables.
Quantitative variables are temporarily omitted to simplify the presentation.
The surrogate models $\tilde{f}$ and $\tilde{g}$, as well as the kernel, are defined on the categorical set $\mathcal{X}^{\cat}$.
A neighborhood constructed at an incumbent $\uu \in \mathcal{X}^{\cat}$ is denoted $\mathcal{N}(\uu;m)$ where $m \in \{1,2,\ldots,|\mathcal{X}^{\cat}|\}$ is the number of neighbors.
This notation is formalized at the end of the section.

\Cref{sec:categorical_distance_objective} presents categorical pseudo-distances that establish proximity with respect to the objective function.
Afterwards, \cref{sec:categorical_distance_constraints} proposes a categorical pseudo-distance for the constraint functions.
Based on these categorical pseudo-distances, surrogate-based neighborhoods are developed in \cref{sec:surrogate_neighborhoods}. 
%
%


\subsection{Categorical pseudo-distances for the objective function}
\label{sec:categorical_distance_objective}

Kernels 
are equivalent to scalar products in Hilbert spaces~\cite{Me1909}.
This property allows
categorical kernels to be transformed into categorical pseudo-distances.
Let $\bm{\phi}:\mathcal{X}^{\cat} \to \mathcal{H}$ be a mapping that maps the categorical components of the categorical set $\mathcal{X}^{\cat}$ into a Hilbert space $\mathcal{H}$ equipped with a scalar product.
Then, the categorical kernel is uniquely defined
according to the Moore-Aronszajn Theorem~\cite{Ar1950} and satisfies $\kappa^{\cat}(\uu,\vv) =  \bm{\phi}(\uu)^{\top} \bm{\phi}(\vv)$.
Note that the categorical kernel $\kappa^{\cat}$ is still computed as described in \cref{sec:kernel}, \textit{i.e.}, with Gaussian kernels or matrices.
{\blr
%
The categorical pseudo-distance $d_f$ is defined from the categorical kernel~\cite{GrBoRaScSm2012} as follows
\begin{align}
\begin{split}
    d_{f}(\uu,\vv) &\coloneq \sqrt{\kappa^{\cat}(\uu,\uu) + \kappa^{\cat}(\vv,\vv) - 2\kappa^{\cat}(\uu,\vv)} \\
    &=  \sqrt{\bm{\phi}(\uu)^{\top} \bm{\phi}(\uu)  + \bm{\phi}(\vv)^{\top} \bm{\phi}(\vv) - 2 \bm{\phi}(\uu)^{\top} \bm{\phi}(\vv)}  \\ 
    &= \| \bm{\phi}(\uu) - \bm{\phi}(\vv) \|_2.
\label{eq:categorical_distance}
\end{split}
\end{align}
{\blr The mapping $\bm{\phi}$} is not necessarily injective on $\mathcal{X}^{\cat}$.
Therefore, two distinct categorical components may have identical feature representations and a zero pseudo-distance.
The function $d_f$ is thus a pseudo-distance in general.
{\blr If $\bm{\phi}$} is injective on $\mathcal{X}^{\cat}$, then $d_f$ is a distance.

}

Although the {\blr mapping $\bm{\phi}$} and the Hilbert space $\mathcal{H}$ are not explicitly known,~\eqref{eq:categorical_distance} shows that the categorical pseudo-distance can be equivalently interpreted in a Hilbert space with $d_{f}(\uu, \vv)=\| \bm{\phi}(\uu) - \bm{\phi}(\vv) \|_2$.
Without constraint functions, categorical neighborhoods are characterized only by this categorical pseudo-distance.
Hence, unconstrained neighborhoods can be viewed as being implicitly constructed in a Hilbert space, where two categorical components $\uu, \vv \in \mathcal{X}^{\cat}$ that are more similar correspond to a smaller pseudo-distance $\| \bm{\phi}(\uu) - \bm{\phi}(\vv) \|_2$.
\Cref{fig:distance_hilbert} illustrates this Hilbert space as a two-dimensional continuous space, where categorical components are mapped into continuous vectors.
Note that this is only an equivalent representation.
{\blr In practice, the pseudo-distance is computed as $d_{f}(\uu,\vv) = \sqrt{\kappa^{\cat}(\uu,\uu) + \kappa^{\cat}(\vv,\vv) - 2\kappa^{\cat}(\uu,\vv)}$.}
%

\begin{figure}[htb!]
\centering
\scalebox{1}{\begin{tikzpicture}[scale=1.2, >=Stealth]

\node at (-3.5,0.75) {$\mathcal{X}^{\cat} = \{ \uu, \vv, \bm{w} \}$};

\draw[->, thick] (-2,0.75) -- (-0.2,0.75) node[midway, above] {\large $\phi$};

\begin{scope}[shift={(0.4,-0.2)}]
  \draw[fill=gray!10, thick] (0,0) rectangle (3,2);
  \node[anchor=west] at (0.05,1.8) {$\mathcal{H}$};

  \coordinate (pu) at (1.5,0.8);
  \coordinate (pv) at (1.675,0.55);
  \coordinate (pw) at (1.8,1.7);

  \fill (pu) circle (2pt);
  \fill (pv) circle (2pt);
  \fill (pw) circle (2pt);

  \node[anchor=south, yshift=2pt] at (pu) {\small$\phi(\uu)$};  
  \node[anchor=west, xshift=2pt] at (pv) {\small$\phi(\vv)$};
  \node[anchor=west, xshift=2pt, yshift=1pt] at (pw) {\small$\phi(\bm{w})$};

\end{scope}

\node[anchor=west] at (3.8,1.15) {$\mathcal{N}(\uu;2)=\mathcal{N}(\vv;2)=\{\uu,\vv\}$};
\node[anchor=west] at (3.8,0.25) {$\mathcal{N}(\bm{w};2)=\{\bm{w},\uu\}$};

\node at (3.8,-0.35) {};

\end{tikzpicture}}
%
\caption{Visualization of {\blr an} unconstrained neighborhood equivalently constructed in a Hilbert space.}
\label{fig:distance_hilbert}
\end{figure}

In \cref{fig:distance_hilbert}, the categorical components are mapped in a Hilbert space via an implicit {\blr mapping $\bm{\phi}$}.
Schematically, neighborhoods containing two elements can be viewed as follows: the neighborhoods constructed at $\uu$ and $\vv$ coincide and contain the components $\{\uu,\vv\}$, while the neighborhood constructed at $\bm{w}$ contains $\{\bm{w},\uu\}$.

On the one hand, the construction of the categorical pseudo-distance $d_{f}$ requires non-negligible computational effort, notably because of the maximization of the marginal likelihood of the GP in~\eqref{eq:log_likelihood}.
On the other hand, it is induced by a finely-tuned kernel.
The categorical pseudo-distance in~\eqref{eq:categorical_distance} is constructed from the problem data through the kernel, so that correlated categorical components are embedded closer to each other in a Hilbert space, while uncorrelated components are placed farther apart.
This behavior is particularly useful for constructing categorical neighborhoods.
This pseudo-distance is finer than binary categorical distances, such as the Gower (or Hamming) distance, which treats all mismatches equally and does not use information from the problem.
As an example, let $\mathcal{X}^{\cat}= \{\text{\sc{r}}, \text{\sc{b}}, \text{\sc{g}} \}$ and $d_{\mathrm{G}}$ denote the Gower distance.
Then $d_{\mathrm{G}}(\text{\sc{r}}, \text{\sc{g}}) = d_{\mathrm{G}}(\text{\sc{r}}, \text{\sc{b}}) = 1$, implying that it is not possible to determine whether green or blue is closer to red.

If a surrogate model is available, but no kernel is defined,
the categorical pseudo-distance can be approximated by another pseudo-distance using prediction values, such that 
$d_{f} \left( \uu, \vv \right) \approx \left| \hat{f}(\uu) -\hat{f}(\vv) \right|$.
This pseudo-distance trivially satisfies symmetry and {\blr triangle} inequality properties, but not the identity of indiscernibles property since $ d_{f} \left( \uu, \vv \right)=0 \not\Leftrightarrow \uu = \vv$ when $\hat{f}$ is not injective.
Hence, categorical components with similar prediction values are considered close.
Conceptually, this mimics the behavior of a kernel inducing an interpolation model, such as GPs. 
Although simple, this pseudo-distance captures proximity with respect to the objective function $f$.
That said, when a kernel with adjusted hyperparameters is available, the categorical pseudo-distance of~\eqref{eq:categorical_distance} is generally preferable, as it corresponds to measuring similarity in an inner-product Hilbert space.

\subsection{Categorical pseudo-distance for constraint functions}
\label{sec:categorical_distance_constraints}

%
As discussed previously, the main goal is to construct neighborhoods that incorporate information from both the objective and constraint functions.
This section develops a pseudo-distance that establishes proximity with respect to the constraint functions.

When a neighborhood is constructed at a feasible component, it should first maintain feasibility and then promote improvement in the objective function, \textit{i.e.}, similarity with respect to the objective.
In this setting, differences in the degree of feasibility between categorical components that are both predicted to be feasible are of limited importance.
From an optimization perspective, once feasibility is predicted, the focus should instead be on improving the objective function.

In particular, a large pseudo-distance between a “highly feasible” component and a “marginally feasible” component should be avoided, as both satisfy the constraints and should be compared primarily regarding the objective.
For this reason, the proposed pseudo-distance is constructed so that it evaluates to zero when both compared components are predicted to be feasible.

The next functions are used to construct the pseudo-distance.
For each $j \in J$, a function $\hat{g}_j^{+} : \mathcal{X}^{\cat} \to \mathbb{R}$ is defined from the predictive model $\hat{g}_j : \mathcal{X}^{\cat} \to \mathbb{R}$ and its associated uncertainty model $\hat{\sigma}_j : \mathcal{X}^{\cat} \to \mathbb{R}_+$, both derived from the surrogate model $\tilde{g}_j$ of the $j$-th constraint.
Let $\lambda \ge 0$ be a relaxation parameter used to relax a predicted constraint
value $\hat{g}_j(\uu)$ using its associated uncertainty
$\hat{\sigma}_j(\uu) \ge 0$.
The function $\hat{g}_j^{+}$ sets to zero any predicted constraint value
$\hat{g}_j(\uu)$ that is feasible after relaxation by
$\lambda \hat{\sigma}_j(\uu)$ and normalizes the remaining values, such that
\begin{align}
    \hat{g}_j^+ \left( \uu \right) \coloneq 
    \begin{cases}
      0 \quad &\text{ if } \hat{g}_j(\uu) - \lambda \hat{\sigma}_j(\uu)   \leq 0, \\
    \psi \left( \hat{g}_j(\uu) \right) &\text{ otherwise}, 
    \end{cases}
\label{eq:surrogate_constraint_plus}
\end{align}
where $\psi$ is a normalization function ensuring that
$\hat{g}_j^+ \left( \uu \right) \in [0,1]$ holds for any
$\uu \in \mathcal{X}^{\cat}$.
For example, $\psi$ may apply a min-max normalization on the bounds of
$\hat{g}_j(\uu)$ when
$\hat{g}_j(\uu) - \lambda \hat{\sigma}_j(\uu) > 0$ is considered.
{\blr
In a blackbox setting, the range of the constraint values is generally unknown.
In practice, an approximate normalization can be obtained from the constraint values available in the current dataset.
For instance, zero can be used as the lower normalization bound and the largest known positive constraint value as the upper normalization bound.
}

Note that if the surrogate model does not provide an uncertainty measure, then
$\lambda$ can be set to zero, meaning that no uncertainty is taken into account.
For readability,~\eqref{eq:surrogate_constraint_plus} may be expressed in a vector form {\blr $\hat{\bm{g}}^+:\mathcal{X}^{\cat} \to [0,1]^{|J|}$}, considering all indices $j \in J$ compactly.

The pseudo-distance function associated with the constraint functions is
constructed from the vector function $\hat{\bm{g}}^+$, as follows: %
\begin{equation*}
    d_{\bm{g}} \left(\uu, \vv \right) \coloneq \lVert \hat{\bm{g}}^+(\uu) - \hat{\bm{g}}^+(\vv) \rVert_p,
\end{equation*}
%
%
{\blr 
where $p \geq 1$.
By construction, $d_{\bm{g}}(\uu,\vv)=0$ whenever both arguments are predicted to be feasible after relaxation, \textit{i.e.}, when $\hat{\bm{g}}(\uu)-\lambda \hat{\boldsymbol{\sigma}}(\uu) \leq \bm{0}$ and $\hat{\bm{g}}(\vv)-\lambda \hat{\boldsymbol{\sigma}}(\vv) \leq \bm{0}$.
}
The normalization function $\psi$ ensures that the pseudo-distance $d_{\bm{g}}$ is not biased by the relative
scales of the constraint function values.
The symmetry and {\blr triangle} inequality properties are trivially satisfied by virtue of the $p$-norm. 
However, the identity of indiscernibles is not respected, \textit{i.e.}, $d_{\bm{g}} \left(\uu, \vv \right) = 0 \not \Leftrightarrow \uu = \vv$ since the vector function $\hat{\bm{g}}^+$ is not injective by design. 
This last unsatisfied property makes $d_{\bm{g}}$ a pseudo-distance function instead of a distance function.
%

\subsection{Surrogate-based neighborhoods}
\label{sec:surrogate_neighborhoods}

The novel categorical neighborhoods are developed in this section.
A neighborhood is constructed at a component $\uu \in \mathcal{X}^{\cat}$.
The other components in the neighborhood are selected based on trade-offs between
pseudo-distances concerning the objective and the constraint functions.
This selection is formalized using notions of dominance inspired by
bi-objective optimization.
In the following, two ranking functions are derived from the pseudo-distance functions introduced in previous sections: a primary and a secondary ranking function, both
parameterized by a component $\uu \in \mathcal{X}^{\cat}$ at which the neighborhood is constructed.
Conceptually, the first selected components correspond to Pareto-optimal
solutions concerning the two ranking functions, while subsequent components
prioritize the primary ranking function over the secondary one.

The primary and secondary ranking functions are defined in two cases depending on the component $\uu \in \mathcal{X}^{\cat}$.
If $\uu$ is feasible, then the primary ranking function is derived
from $d_{\bm{g}}$ (constraints),
since the priority is to construct neighborhoods with components that maintain feasibility.
Otherwise, $\uu$ is infeasible, and the primary ranking function is based on $d_f$ (objective), allowing greater flexibility with respect to feasibility.
This may allow finding promising solutions
that are difficult to reach when imposing feasibility at all times.

%
The primary and secondary ranking functions, noted $p_{\uu}:\mathcal{X}^{\cat} \to \mathbb{R}_+$ and $s_{\uu}: \mathcal{X}^{\cat} \to \mathbb{R}_+$ respectively, are parameterized by $\uu \in \mathcal{X}^{\cat}$ and defined as follows
\begin{align*}
    \begin{split}
        &p_{\uu}( \vv) \coloneq d_{\bm{g}} \left( \uu, \vv \right) \ \text{ and } \ s_{\uu}( \vv) \coloneq d_{f} \left( \uu, \vv \right) \ \text{ if } \uu \in \Omega, \\
        &p_{\uu}( \vv) \coloneq d_{f} \left( \uu, \vv \right) \ \text{ and } \ s_{\uu}( \vv) \coloneq d_{\bm{g}} \left( \uu, \vv \right) \ \text{ if } \uu \not \in \Omega. \\
    \end{split}
\end{align*}
%
%

%
%
%
%
{\blr
An illustrative example is presented before formally defining the ordering relation used to construct a neighborhood.
In this example, there are twelve categorical components.
{\blr 
The neighborhood is constructed around a component $\uu \in \mathcal{X}^{\cat}$, and the other eleven components are ranked.
}
For a constrained problem, the ranking is performed in three successive steps, as illustrated from left to right in \cref{fig:rank}.
}

\begin{figure}[htb!]
\centering
%
\begin{minipage}[b]{0.75\textwidth}
  \centering

  \begin{subfigure}[t]{0.32\textwidth}
    \centering
    \captionsetup{justification=centering}
    \scalebox{0.9}{\begin{tikzpicture}[baseline=(current bounding box.south), x=0.58cm,y=0.48cm,>=Latex, every node/.style={font=\normalsize}]
    \draw[-{Latex}] (0,0) -- (6.25,0) node[below left=1pt] {$p_{\uu}$};
  \draw[-{Latex}] (0,0) -- (0,6) node[left=1pt] {$s_{\uu}$};

  \fill (0, 2) circle (1.6pt) node[left=1pt] {1};
  \fill (1, 1.5) circle (1.6pt) node[left=1pt] {2};
  \fill (2, 1) circle (1.6pt) node[left=1pt] {3};
  \fill (3, 0.5) circle (1.6pt) node[left=1pt] {4};

  \fill (0,4) circle (1.6pt);
  \fill (0,5) circle (1.6pt);
  \fill (1.90,2.55) circle (1.6pt);
  \fill (2.70,3.65) circle (1.6pt);
  \fill (3.30,3.20) circle (1.6pt);
  \fill (4.00,2.25) circle (1.6pt);
  \fill (5.00,3.45) circle (1.6pt);
\end{tikzpicture}}
    \caption{Pareto components.}
    \label{subfig:rank1}
  \end{subfigure}%
  \hfill
  \begin{subfigure}[t]{0.32\textwidth}
    \centering
    \captionsetup{justification=centering}
    \scalebox{0.9}{\begin{tikzpicture}[baseline=(current bounding box.south), x=0.58cm,y=0.48cm,>=Latex, every node/.style={font=\normalsize}]
    \draw[-{Latex}] (0,0) -- (6.25,0) node[below left=1pt] {$p_{\uu}$};
  \draw[-{Latex}] (0,0) -- (0,6) node[left=1pt] {$s_{\uu}$};

  \fill (0,4) circle (1.6pt) node[left=1pt] {5};
  \fill (0,5) circle (1.6pt) node[left=1pt] {6};

  \fill (0, 2) circle (1.6pt) node[left=1pt] {1};
  \fill (1, 1.5) circle (1.6pt) node[left=1pt] {2};
  \fill (2, 1) circle (1.6pt) node[left=1pt] {3};
  \fill (3, 0.5) circle (1.6pt) node[left=1pt] {4};

  \fill (1.90,2.55) circle (1.6pt);
  \fill (2.70,3.65) circle (1.6pt);
  \fill (3.30,3.20) circle (1.6pt);
  \fill (4.00,2.25) circle (1.6pt);
  \fill (5.00,3.45) circle (1.6pt);
\end{tikzpicture}}
    \caption{Non-Pareto components with $p_{\uu}(\vv)=0$.}
    \label{subfig:rank2}
  \end{subfigure}%
  \hfill
  \begin{subfigure}[t]{0.32\textwidth}
    \centering
    \captionsetup{justification=centering}
    \scalebox{0.9}{\begin{tikzpicture}[baseline=(current bounding box.south), x=0.58cm,y=0.48cm,>=Latex, every node/.style={font=\normalsize}]
    \draw[-{Latex}] (0,0) -- (6.25,0) node[below left=1pt] {$p_{\uu}$};
  \draw[-{Latex}] (0,0) -- (0,6) node[left=1pt] {$s_{\uu}$};

  \fill (0,4) circle (1.6pt) node[left=3pt] {5};
  \fill (0,5) circle (1.6pt) node[left=3pt] {6};

  \fill (0, 2) circle (1.6pt) node[left=1pt] {1};
  \fill (1, 1.5) circle (1.6pt) node[left=1pt] {2};
  \fill (2, 1) circle (1.6pt) node[left=1pt] {3};
  \fill (3, 0.5) circle (1.6pt) node[left=1pt] {4};

  \fill (1.90,2.55) circle (1.6pt) node[left=1pt] {7};
  \fill (2.70,3.65) circle (1.6pt) node[left=1pt] {8};
  \fill (3.30,3.20) circle (1.6pt) node[left=1pt] {9};
  \fill (4.00,2.25) circle (1.6pt) node[left=1pt] {10};
  \fill (5.00,3.45) circle (1.6pt) node[left=1pt] {11};
\end{tikzpicture}}
    \caption{Remaining components ordered with the primary ranking function.}
    \label{subfig:rank3}
  \end{subfigure}

\end{minipage}%
\hspace{1.5em}
\tikz[baseline]{
  \draw[draw=black,
        line width=0.8pt,
        dash pattern=on 0pt off 2.2pt,
        line cap=round]
        (0,-1cm) -- (0,3.6cm);
}
\hspace{0.2em}
%
\begin{minipage}[b]{0.15\textwidth}
  \centering

  \begin{subfigure}[t]{\textwidth}
    \centering
    \captionsetup{justification=centering}
    \scalebox{0.9}{\begin{tikzpicture}[baseline=(current bounding box.south), x=0.58cm,y=0.48cm,>=Latex, every node/.style={font=\normalsize}]
  \draw[-{Latex}] (0,0) -- (3.5,0) node[below left=1pt] {$p_{\uu}$};
  \draw[-{Latex}] (0,0) -- (0,6) node[above left=2pt] {$s_{\uu}$};

  \fill (0,1)    circle (1.6pt) node[left=3pt] {1};
  \fill (0,2.5)  circle (1.6pt) node[left=3pt] {2};
  \fill (0,4)    circle (1.6pt) node[left=3pt] {3};
  \fill (0,4.75) circle (1.6pt) node[left=3pt] {4};
  \fill (0,5.75) circle (1.6pt) node[left=3pt] {5};
\end{tikzpicture}}
    \caption{Unconstrained problem.}
    \label{subfig:rank_unconstrained}
  \end{subfigure}

\end{minipage}

\caption{Ordering components with ranking functions. (\subref{subfig:rank1})-(\subref{subfig:rank3}) three steps for constrained problems. (\subref{subfig:rank_unconstrained}) single step for unconstrained problems.
The component $\uu$ is located at the origin and is not displayed.}
\label{fig:rank}
\end{figure}

The first four components are the Pareto components, \textit{i.e.}, the nondominated components, in the space of images.
The Pareto components are ordered among themselves solely with the primary ranking function $p_{\uu}$.
The components~5 and~6 are not Pareto, but they have a primary ranking function value of zero with respect to the component $\uu$ at which the neighborhood is constructed.
This situation may most notably occur when the primary ranking function is based on the pseudo-distance $d_{\bm{g}}$: recall that when both components are predicted to be feasible, this pseudo-distance returns zero. 
Among themselves, components~5 and~6 are ordered using the secondary ranking function.
Finally, the remaining components~7 to~11 are directly ordered with the primary ranking function.

Note that when the problem is unconstrained, the component $\uu \in \Omega$ and, consequently, the primary ranking function is based on $d_{\bm{g}}$.
However, in this case, the pseudo-distance function $d_{\bm{g}}$ is always set to zero by convention, since its components are always feasible.
Therefore, only the pseudo-distance function of the objective $d_f$ is considered.
An unconstrained example is illustrated with six components (including $\uu$ that is not displayed) in \cref{subfig:rank_unconstrained}, where the ordering is done directly based on their height.

{\blr
The discussion above suggests grouping the categorical components into three disjoint subsets, each illustrated by the group of numbered points shown in one of
\cref{subfig:rank1,subfig:rank2,subfig:rank3}.
The subsets are considered from left to right, which directly determines their relative ordering.
For example, any component represented in \cref{subfig:rank1} is ordered before any component represented in \cref{subfig:rank2}.
Within the subsets shown in \cref{subfig:rank1,subfig:rank3}, the components are ordered using the primary ranking function $p_{\uu}$, from left to right.
Within the subset shown in \cref{subfig:rank2}, the components are ordered using the secondary ranking function $s_{\uu}$, from bottom to top.
The following definition formalizes this ordering.
}

%
%
%
%
    %




    
%
%
%
{\blr
%
\begin{definition}[Ordering relation]
For a given categorical component $\uu \in \mathcal{X}^{\cat}$, let
$\mathcal{C}_{\uu} \coloneq \mathcal{X}^{\cat}\setminus\{\uu\}$
denote the categorical components to be ranked, and let
$p_{\uu}:\mathcal{X}^{\cat}\to\mathbb{R}_+$
and
$s_{\uu}:\mathcal{X}^{\cat}\to\mathbb{R}_+$
denote the associated primary and secondary ranking functions.
%

The following three disjoint subsets are defined, where the first subset contains the components that are nondominated with respect to $(p_{\uu},s_{\uu})$:
\begin{align*}
\mathcal{A}_{\uu}^{1}
&\coloneq
\left\{
\vv\in\mathcal{C}_{\uu} :
\nexists\,\w\in\mathcal{C}_{\uu}
\text{ such that }
p_{\uu}(\w)\leq p_{\uu}(\vv),
\;
s_{\uu}(\w)\leq s_{\uu}(\vv),
\text{ with at least one strict inequality}
\right\},
\\
\mathcal{A}_{\uu}^{2}
&\coloneq
\left\{
\vv\in\mathcal{C}_{\uu}\setminus\mathcal{A}_{\uu}^{1} :
p_{\uu}(\vv)=0
\right\},
\\
\mathcal{A}_{\uu}^{3}
&\coloneq
\mathcal{C}_{\uu}
\setminus
\left(
\mathcal{A}_{\uu}^{1}
\cup
\mathcal{A}_{\uu}^{2}
\right).
\end{align*}
%
%
%
%
%
%
%
%
%
The ordering relation is such that, for any pair
$\vv,\w\in\mathcal{C}_{\uu}$, $\vv$ is ordered before or tied to $\w$,
denoted $\vv\preceq_{\uu}\w$, when one of the following mutually exclusive conditions holds:
\begin{enumerate}
    \item $\vv\in\mathcal{A}_{\uu}^{i}$ and
    $\w\in\mathcal{A}_{\uu}^{j}$, with
    $i,j\in\{1,2,3\}$ such that $i<j$;

    \item $\vv,\w\in\mathcal{A}_{\uu}^{i}$, with
    $i\in\{1,3\}$, such that
    $p_{\uu}(\vv)\leq p_{\uu}(\w)$;

    \item $\vv,\w\in\mathcal{A}_{\uu}^{2}$, such that
    $s_{\uu}(\vv)\leq s_{\uu}(\w)$.
\end{enumerate}

\label{def:ordering_relation}
\end{definition}
}

{\blr
The ordering relation in \cref{def:ordering_relation} defines a total preorder, but not necessarily a total order, since ties may occur between components within the same subset.
Indeed, the three subsets form a partition of $\mathcal{C}_{\uu}$ and are considered successively, while the ordering within each subset is induced by the usual order on $\mathbb{R}$.
Consequently, the relation is reflexive, transitive, and total.
Nevertheless, it provides a systematic way to identify promising categorical components using information from both the objective and constraint functions.
}

%

Now that the ordering relation is formally defined, the neighborhoods based on surrogate models are finally introduced.

{\blr
\begin{definition}[Surrogate-based neighborhood]
For a given component $\uu \in \mathcal{X}^{\cat}$ and integer
$m \in \{1,2,\ldots,|\mathcal{X}^{\cat}|\}$, the
{\em surrogate-based neighborhood}
$\mathcal{N}(\uu;m) \subseteq \mathcal{X}^{\cat}$
contains $\uu$ and the $m-1$ components of $\mathcal{C}_{\uu}$ with lowest ordering according to $\preceq_{\uu}$.
The relation $\vv \preceq_{\uu} \w$ is satisfied whenever
$\vv \in \mathcal{N}(\uu;m)\setminus\{\uu\}$ and
$\w \notin \mathcal{N}(\uu;m)$. 
If a tie occurs, the component encountered first when enumerating
$\mathcal{X}^{\cat}$ is selected.
%
\label{def:neighborhood}
\end{definition}
}

{\blr
The enumeration rule is only used to break ties and ensure that the neighborhood is uniquely defined.
Other deterministic tie-breaking rules could be used without changing the ordering relation $\preceq_{\uu}$.
}

{\blr 
The current neighborhood construction enumerates all
$\left|\mathcal{X}^{\cat}\right|$ categorical components.
For constrained problems, identifying the nondominated components requires
$\mathcal{O}\left(\left|\mathcal{X}^{\cat}\right|^2\right)$
pairwise comparisons.
The subsequent ordering requires at most
$\mathcal{O}\left(\left|\mathcal{X}^{\cat}\right|
\log \left|\mathcal{X}^{\cat}\right|\right)$
operations, which essentially consists of sorting the components within their respective subsets.
For unconstrained problems, the components are ordered directly according to
$d_f$, requiring
$\mathcal{O}\left(\left|\mathcal{X}^{\cat}\right|
\log \left|\mathcal{X}^{\cat}\right|\right)$
operations.
Hence, the overall worst-case complexity of constructing a neighborhood is
$\mathcal{O}\left(\left|\mathcal{X}^{\cat}\right|^2\right)$.
}

\Cref{def:neighborhood} generalizes distance-induced neighborhoods
in~\cite{catmads} by using relations instead of a single distance for selecting
components.
In fact, when a problem is unconstrained, a surrogate-based neighborhood is equivalent to a distance-induced neighborhood using $d_{f}$.
That said, the implementation in~\cite{catmads} uses a basic mixed-variable interpolation to construct the categorical distance.
{\blr
%
The kernel-induced categorical pseudo-distance $d_f$ introduced in \cref{sec:categorical_distance_objective}  provides a more structured similarity measure, as it is induced by an inner product in a Hilbert space through well-defined categorical kernels.
The proposed improvements are therefore not limited to constrained problems, since the kernel-induced pseudo-distance also compares categorical components through the geometry induced by the kernel in the unconstrained case.
}

\subsection{Illustrative case study of the mechanical-part problem}
\label{sec:illustrative_example}

The mechanical-part design problem from \cref{sec:example} is used here to develop a case study with categorical neighborhoods.
{\blr
A single Latin Hypercube Sampling (LHS) of 96 points is performed on the mixed-variable domain described in \cref{sec:example},
 in which four points are generated for each of the 24 categorical components.
}
Only one sampled point is feasible, that is, $\x_{\star} = \left(\B, \Wood, \Circle, 5.5, 1 \right)$, with objective function value $f(\x_{\star})=28.55$.
{\blr
In this study, the continuous variables are fixed at their values in $\x_{\star}$.
The goal is to identify another categorical component that remains feasible and improves the objective value.
}
The continuous variables are fixed, and the current categorical component is noted $\uu = \left(\B, \Wood, \Circle\right) \in \mathcal{X}^{\cat}$.

{\blr
For two categorical components $\uu,\vv\in\mathcal{X}^{\cat}$, the Gower distance is defined as
\begin{equation*}
d_{\mathrm{G}}(\uu,\vv)
=
\sum_{i=1}^{n^{\cat}}
\mathbb{I}\left[u_i\neq v_i\right],
\end{equation*}
where $\mathbb{I}[\cdot]$ is the indicator function.
It counts the number of categorical variables that differ between two components, \textit{e.g.}, $
d_{\mathrm{G}}\left(
(\B,\Wood,\Circle),
(\A,\Wood,\Circle)
\right)=1$.
The neighborhood induced by this distance contains the six categorical components at distance one from
$\uu=(\B,\Wood,\Circle)$, \textit{i.e.}, the components obtained by changing exactly one categorical variable.
}
\Cref{tab:neighbors_comparison} compares three neighborhood-generation strategies
    producing exactly six points.
{\blr Each point has the same continuous variables as $\x_{\star}$ and is different from $\uu$.}
The first and simplest strategy uses the Gower distance \cite{Gower1971}. 
%
The second strategy uses an objective-based categorical distance~\cite{catmads}, 
    induced by a mixed-variable Inverse Distance Weighting (IDW) regression model.
The third strategy employs the surrogate-based neighborhoods from \cref{def:neighborhood}, using GPs as surrogate models.

\begin{table}[htb!]
\blr
\centering
\begin{tabular}{c c c c c c c c c}
\hline
Strategy & Neighbor & Supplier & Material & Shape & $f$ & $g_1$ & $g_2$ & Feasible \\
\hline
\multirow{6}{*}{Gower distance~\cite{Gower1971}}
& 1 & $\B$ & $\Aluminum$ & $\Circle$ & 18.4928 & -0.1300 & 0.0150 & No \\
& 2 & $\B$ & $\Steel$ & $\Circle$ & 14.5128 & -0.0225 & 0.0625 & No \\
& 3 & $\B$ & $\Composite$ & $\Circle$ & 9.5329 & -0.0800 & 0.0575 & No \\
& 4 & $\B$ & $\Wood$ & $\Square$ & 24.0428 & 0.0050 & 0.0650 & No \\
& 5 & $\B$ & $\Wood$ & $\Ellipse$ & 24.5628 & -0.1075 & 0.0125 & No \\
& 6 & $\A$ & $\Wood$ & $\Circle$ & 28.3468 & 0.3450 & 0.5600 & No \\
\hline
\multirow{6}{*}{Objective-based~\cite{catmads}}
& 1 & $\B$ & $\Wood$ & $\Square$ & 24.0428 & 0.0050 & 0.0650 & No \\
& 2 & $\A$ & $\Wood$ & $\Circle$ & 28.3468 & 0.3450 & 0.5600 & No \\
& 3 & $\B$ & $\Steel$ & $\Circle$ & 14.5128 & -0.0225 & 0.0625 & No \\
& 4 & $\B$ & $\Aluminum$ & $\Circle$ & 18.4928 & -0.1300 & 0.0150 & No \\
& 5 & $\B$ & $\Wood$ & $\Ellipse$ & 24.5628 & -0.1075 & 0.0125 & No \\
& 6 & $\B$ & $\Composite$ & $\Circle$ & 9.5329 & -0.0800 & 0.0575 & No \\
\hline
\multirow{6}{*}{Surrogate-based}
& $\bm{1}$ & $\text{\bfseries\scshape b}$ & $\text{\bfseries\scshape steel}$ & $\text{\bfseries\scshape square}$ & $\bm{10.0028}$ & $\bm{-0.0225}$ & $\bm{-0.0451}$ & \textbf{Yes} \\
& 2 & $\B$ & $\Wood$ & $\Ellipse$ & 24.5628 & -0.1075 & 0.0125 & No \\
& 3 & $\B$ & $\Aluminum$ & $\Circle$ & 18.4928 & -0.1300 & 0.0150 & No \\
& $\bm{4}$ & $\text{\bfseries\scshape a}$ & $\text{\bfseries\scshape alum}$ & $\text{\bfseries\scshape circle}$ & $\bm{18.2868}$ & $\bm{-0.0225}$ & $\bm{-0.0451}$ & \textbf{Yes} \\
& $\bm{5}$ & $\text{\bfseries\scshape a}$ & $\text{\bfseries\scshape wood}$ & $\text{\bfseries\scshape ellipse}$ & $\bm{24.3568}$ & $\bm{-0.6900}$ & $\bm{-0.4715}$ & \textbf{Yes} \\
& $\bm{6}$ & $\text{\bfseries\scshape a}$ & $\text{\bfseries\scshape comp}$ & $\text{\bfseries\scshape square}$ & $\bm{8.8169}$ & $\bm{-0.0225}$ & $\bm{-0.0451}$ & \textbf{Yes} \\
\hline
\end{tabular}
\caption{Neighbors of $\uu=\left(\B, \Wood, \Circle\right)$ proposed by three strategies.}
\label{tab:neighbors_comparison}
\end{table}

The results highlight clear differences between the three strategies.
The Gower distance and the objective-based neighborhood fail to identify any feasible point, whereas the surrogate-based neighborhood proposed identifies four feasible points.
Recall that \cref{subfig:scatter3} shows the results of a Latin Hypercube Sampling with $100{,}000$ samples per categorical component, suggesting that only six categorical components are likely to admit feasible continuous regions.
The proposed method identifies four of these six components within the six evaluations, while the incumbent component accounts for another one.
The only remaining feasible component is $\left(\B, \Aluminum, \Square\right)$, which appears to be extremely unlikely to yield feasible points, with only $0.01\%$ feasibility observed in {\blr an} LHS of $100{,}000$ points in \cref{subfig:scatter3}.

Only the surrogate-based strategy improves the objective value while maintaining feasibility.
The best improvement is obtained with the categorical component $\left(\A, \Composite, \Square\right)$, yielding an objective value of $8.8169$.
This component also corresponds to that of the best solution identified in \cref{sec:example} using an LHS with $100{,}000$ points per component.

The ordering and selection of categorical components in the surrogate-based neighborhood are presented in \cref{fig:neighbors_example}.
The figure plots the primary and secondary ranking functions $p_{\uu}$ and~$s_{\uu}$ introduced in \cref{sec:surrogate_neighborhoods}, with logarithmic scales for readability.
The selected components are shown in blue, together with their categorical component and a number indicating their ranking.
On the horizontal axis $p_{\uu}$, an axis break is introduced to display components with distance zero on a logarithmic scale.
In the figure, the feasible neighbors all have a distance of zero with respect to $p_{\uu}$: they lie on the vertical axis.
%
%
{\blr The surrogate-based approach ranks $\left(\B, \Steel, \Square\right)$ first, whereas this component is not selected by the objective-based neighborhood.}

\begin{figure}[htb!]
\centering
\def\PUZERO{1e-5}

\begin{tikzpicture}
\begin{loglogaxis}[
    width=1\linewidth,
    height=9cm,
    xmin=\PUZERO,
    xmax=1,
    ymin=1e-6,
    ymax=1e-1,
    grid=major,
    major grid style={gray!35},
    minor tick num=0,
    tick align=outside,
    clip=false,
    xtick={1e-4,1e-3,1e-2,1e-1,1},
    extra x ticks={\PUZERO},
    extra x tick labels={$0$},
    extra x tick style={
        major tick length=0pt,
        tick label style={yshift=2pt}
    },
]

\addplot[
    only marks,
    mark=*,
    mark size=1.8pt,
    color=black
] coordinates {
(5.5024986777147372e-02,7.8089796519540400e-03)
(1.9997524909457470e-02,4.0730123931131956e-05)
(\PUZERO,7.7751486083585952e-03)
(1.4977504335508561e-02,6.7671612646869050e-06)
(1.9996236129080802e-02,7.7829551742105707e-03)
(4.2458889974025991e-02,1.4604167647203070e-05)
(3.0024329308403125e-02,7.8089796519540400e-03)
(\PUZERO,4.0730123931131956e-05)
(6.2571577866207123e-02,7.7751486083585952e-03)
(6.2245669007587517e-02,6.7671612646869050e-06)
(2.9915266726462705e-02,7.7829551742105707e-03)
(3.0013178591780473e-02,1.4604167647203070e-05)
(\PUZERO,3.9814515997971256e-02)
(3.0077836011516097e-02,3.2171067249630036e-02)
(4.2515079871554770e-02,3.9781228466862384e-02)
(5.6337622054846617e-02,3.2137649918795530e-02)
(5.2482573797107614e-02,3.9788909616363943e-02)
(2.2454073420322629e-02,3.2145361019785490e-02)
(7.9999787736330166e-02,7.8022388901937134e-03)
(6.5012280505818906e-02,3.3963077583187840e-05)
(6.5772936168982488e-01,7.7684077321276135e-03)
(\PUZERO,7.7762143243940152e-03)
};

\addplot[
    only marks,
    mark=*,
    mark size=1.8pt,
    color=blue,
    nodes near coords,
    point meta=explicit,
    every node near coord/.append style={font=\large, black, yshift=2pt}
] coordinates {
(1.05e-02,1.2e-05) [2]
(1.4977504335508561e-02,6.7671612646869050e-06) [3]
};

\addplot[
    only marks,
    mark=*,
    mark size=1.8pt,
    color=blue
] coordinates {
(\PUZERO,4.0730123931131956e-05)
(\PUZERO,7.7751486083585952e-03)
(\PUZERO,7.7762143243940152e-03+0.01)
(\PUZERO,3.9814515997971256e-02)
};

\addplot[
    only marks, mark=none, nodes near coords, point meta=explicit,
    every node near coord/.append style={font=\large, black, anchor=east, xshift=0pt, yshift=2pt}
] coordinates {(\PUZERO,4.0730123931131956e-05) [1]};

\addplot[
    only marks, mark=none, nodes near coords, point meta=explicit,
    every node near coord/.append style={font=\large, black, anchor=east, xshift=0pt, yshift=-2pt}
] coordinates {(\PUZERO,7.7751486083585952e-03) [4]};

\addplot[
    only marks, mark=none, nodes near coords, point meta=explicit,
    every node near coord/.append style={font=\large, black, anchor=east, xshift=0pt, yshift=6pt}
] coordinates {(\PUZERO,7.7762143243940152e-03+0.01) [5]};

\addplot[
    only marks, mark=none, nodes near coords, point meta=explicit,
    every node near coord/.append style={font=\large, black, anchor=east, xshift=0pt, yshift=3pt}
] coordinates {(\PUZERO,3.9814515997971256e-02) [6]};

\begin{scope}[xshift=481pt]
    \node[blue, font=\small, xshift=51pt, yshift=0pt]   at (axis cs:1e-11,4.073e-05) {(\B, \Steel, \Square)};
    \node[blue, font=\small, xshift=48pt, yshift=15pt]   at (axis cs:1e-11,7.775e-03) {(\A, \Aluminum, \Circle)};
    \node[blue, font=\small, xshift=50pt, yshift=0pt]  at (axis cs:1e-11,7.776e-03) {(\A, \Wood, \Ellipse)};
    \node[blue, font=\small, xshift=50pt, yshift=0pt]   at (axis cs:1e-11,3.981e-02) {(\A, \Composite, \Square)};
\end{scope}
\begin{scope}[xshift=-6pt, yshift=0pt]
    \node[blue, font=\small, xshift=-45pt, yshift=5pt] at (axis cs:1.05e-02,1.2e-05) {(\B, \Wood, \Ellipse)};
    \node[blue, font=\small, xshift=-45pt, yshift=-2pt] at (axis cs:1.498e-02,6.767e-06) {(\B, \Aluminum, \Circle)};
\end{scope}

\end{loglogaxis}

\draw[black, line width=0.9pt]
  ([xshift=11pt,yshift=-4pt]current axis.south west) --
  ([xshift=19pt,yshift=5pt]current axis.south west);

\draw[black, line width=0.9pt]
  ([xshift=15pt,yshift=-4pt]current axis.south west) --
  ([xshift=23pt,yshift=5pt]current axis.south west);

\draw[-{Latex}, line width=0.9pt]
  (current axis.south west) -- ([xshift=12pt]current axis.south east)
  node[below left=1pt, xshift=16pt] {$p_u$};

\draw[-{Latex}, line width=0.9pt]
  (current axis.south west) -- ([yshift=12pt]current axis.north west)
  node[left=1pt] {$s_u$};

\end{tikzpicture}
\caption{Ordering of categorical components with the surrogate-based neighborhood in the mechanical-part design problem.
The incumbent $\uu=\left(\B, \Wood, \Circle\right)$ is at the origin and it is not displayed.
}
\label{fig:neighbors_example}
\end{figure}

\section{Adapting surrogate-based neighborhoods for mixed-variable direct search}
\label{sec:direct_search}

The categorical neighborhoods developed in \cref{sec:neighborhoods} are now used for optimization in this section, which also reintroduces the integer and continuous variables {\blr ignored} previously.
%
%

\subsection[Adapting the surrogate-neighborhoods for mixed-variable domains]
{Adapting the surrogate-neighborhoods for {\blr mixed-variable} domains}

The pseudo-distances from \cref{sec:neighborhoods} require minor
adaptations now that the surrogate models are defined on the domain
$\mathcal{X}$ rather than on the categorical set $\mathcal{X}^{\cat}$.
They are now computed between points in the domain, while the quantitative
variables are held fixed so that comparisons only concern the categorical
variables.
More precisely, for a given and fixed quantitative component
{\blr $\x^{\quant} \in \mathcal{X}^{\quant}$}, the pseudo-distances introduced
in \cref{sec:neighborhoods} are computed on the set
%
%
\begin{equation*}
(\mathcal{X} \times \mathcal{X})_{\x^{\quant}}
\coloneq
\left\{ (\x,\y) \in \mathcal{X} \times \mathcal{X} :
\x^{\quant} = \y^{\quant} \right\}.
\end{equation*}

The pseudo-distance functions are said to be parameterized by the quantitative component $\x^{\quant}$,
since it influences their behavior without appearing as an explicit argument.

Consequently, the primary and secondary ranking functions $p_{\x}:\mathcal{X} \to \mathbb{R}_+$ and $s_{\x} : \mathcal{X} \to {\blr \mathbb{R}_+}$ are now defined on the domain $\mathcal{X}$ and parameterized at a point $\x \in \mathcal{X}$, which can be either feasible or infeasible.
In this context, an argument of the ranking functions shares the same quantitative component as the incumbent solution, that is $\x^{\quant}$.
The ordering relation $\preceq$ thus compares points of the domain sharing the same quantitative component, making the comparisons effectively on categorical components.

Finally, surrogate-based neighborhoods must be generalized to the mixed-variable
setting.
{\blr
The ordering relation is defined over the full mixed-variable domain, but the neighborhoods are restricted to categorical components.
%
This choice ensures consistency with the \catmads framework.
}
%
%
The following definition adapts \cref{def:neighborhood} in the more general mixed-variable case.

\begin{definition}[Mixed-variable surrogate-based neighborhood]
For a given point $\x \in \mathcal{X}$ and integer $m \in \{1,2,\ldots, \left| \mathcal{X}^{\cat} \right|\}$, the {\em surrogate-based neighborhood} $\mathcal{N}(\x;m) \subseteq \mathcal{X}^{\cat}$ is the set containing $m$ categorical components of
the lowest ordering with respect to an ordering relation $\preceq$ defined on $\mathcal{X}$.
The relation $\y \preceq \z$ is satisfied whenever 
 $$\begin{array}{clcl} &
    \y=\left(\y^{\cat}, \x^{\quant} \right) \in \mathcal{X}
    & \mbox{ where } &
    \y^{\cat} \in \mathcal{X}^{\cat} \cap \mathcal{N}(\x;m)\\
    \mbox{ and } &
    \z=\left(\z^{\cat}, \x^{\quant} \right) \in \mathcal{X}
    & \mbox{ where } &
    \z^{\cat} \in \mathcal{X}^{\cat} \setminus \mathcal{N}(\x;m).
    \end{array}$$
\label{def:neighborhood_general}
\end{definition}
Note that \Cref{def:neighborhood_general} reduces to \cref{def:neighborhood} when restricted to categorical components.
Moreover, a consequence of this definition is that for any $\x \in \mathcal{X}$ and $m \in \{1,2,\ldots, \left| \mathcal{X}^{\cat} \right|\}$,
\begin{equation*}
    \mathcal{N}\left( \x; m \right) = \left\{ \y_{(1)}^{\cat}, \, \y_{(2)}^{\cat}, \ldots, \y_{(m)}^{\cat}  \right\} \subseteq \mathcal{X}^{\cat}
\end{equation*}
in which  $\y_{(1)}^{\cat} = \x^{\cat}$
 and where 
 $\left( \y_{(1)}^{\cat}, \x^{\quant} \right) \preceq 
  \left( \y_{(2)}^{\cat}, \x^{\quant} \right) \preceq \ldots \preceq 
  \left( \y_{(m)}^{\cat}, \x^{\quant} \right)$.
In the next sections, the surrogate-based neighborhoods from \cref{def:neighborhood_general} 
are used to handle the categorical variables.

\subsection{Background on \catmads}
\label{sec:CatMADS}

{\blr 
In this work, \catmads\ is used as an application for the proposed categorical neighborhoods.
Only a brief overview of the method is provided here: see~\cite{catmads} for a complete description.
}
\catmads is a direct search framework that iteratively seeks points yielding a strict decrease of the objective function and/or the constraints~\cite{catmads}.
It generalizes \mads by tackling categorical variables with neighborhoods induced by data from the problem.
The quantitative variables are essentially treated as they are in \mads.
For simplicity, \catmads is presented in \cref{algo:cat_mads_algo} without considering constraints.

\begin{algorithm}[htb!]
\small

0.~\textbf{Initialization}. Set $k=0$, perform {\blr a} Design of Experiment (DoE) and define an  initial mesh \;
\vspace{0.01cm}

1.~\textbf{Opportunistic search} (optional). \;

\vspace{0cm}

2.~\textbf{Opportunistic poll}. Perform polling around the incumbent solution \; 

\vspace{0cm}

3.~\textbf{Opportunistic extended poll} (optional). If Steps 1 \& 2 are unsuccessful, perform quantitative polls  \;
%
\hspace{0.275cm} around points in the categorical poll with objective function values sufficiently close to $f( \x_{(k)} )$ \;

\vspace{0cm}

4.~\textbf{Update}. Set $k\leftarrow k+1$, update mesh and check stopping criterion \;

\hspace{0.275cm}  If the iteration is successful, increase mesh size and \textbf{GO TO} Step 1. \;

\hspace{0.275cm} Else, decrease mesh size, and if the mesh is at its minimum size, then \textbf{STOP}

\caption{The \catmads framework for unconstrained problems (adapted from~\cite{catmads}).}
\label{algo:cat_mads_algo}
\end{algorithm}

Before any optimization is done, a DoE samples points in the domain using a portion of the budget of evaluations.
The DoE is used to gather information on the problem and construct categorical neighborhoods, as well as for determining an initial solution.

The main steps of the \catmads algorithm are Steps 1, 2, and 3. 
If a solution with strictly lower objective function value is found during these steps, it becomes the incumbent solution and the iteration is deemed \textit{successful}.
Otherwise, the iteration is \textit{unsuccessful}.
The \textit{search} is an optional step allowing flexible evaluation of points at the intersection of the mesh and the domain. 
The poll is the core mechanism of \catmads.
At iteration $k$, trial points are generated around the incumbent $\x_{(k)}$ through two polls, a quantitative one and a categorical one.
The quantitative poll fixes the categorical component and performs a standard \mads poll on the quantitative variables.
{\blr
The categorical poll fixes the quantitative component and explores neighboring categorical components selected through a surrogate-based ordering relation.
}
All categorical components in the selected neighborhood, except the incumbent,
form the categorical poll and are evaluated opportunistically.
Two consecutive quantitative polls on $\mathbb{R}^2$ are illustrated in
\cref{subfig:MADS1,subfig:MADS2}.
\Cref{subfig:CatMADS} illustrates both poll types.

\begin{figure}[htb!]
\centering
%
\begin{minipage}[b]{0.68\textwidth}
  \centering
  \begin{subfigure}[t]{0.49\textwidth}
    \centering
    \captionsetup{justification=centering}
    \scalebox{1.05}{\begin{tikzpicture}
        \draw[step=0.5cm,gray,very thin] (0,0) grid (4,4);
        \draw[thick] (1,1) rectangle (3,3);

        \draw[dashed, -{Latex[length=3mm, width=2mm]}, thick] (2,2) -- (1,2);   
        \fill[black] (1,2) circle (2pt);                                        
        \draw[dashed, -{Latex[length=3mm, width=2mm]}, thick] (2,2) -- (2,3);   
        \fill[black] (2,3) circle (2pt);                                        
        \draw[dashed, -{Latex[length=3mm, width=2mm]}, thick] (2,2) -- (2.5,1.5); 
        \fill[black] (2.5,1.5) circle (2pt);                                     

        \fill[black] (2,2) circle (2pt) node[above left] {$\x_{(k)}^{\quant}$};

        \node[left] at (1,2.25) {$\y_{(1)}^{\quant}$};
        \node[below] at (2.7,1.65) {$\y_{(2)}^{\quant}$};
        \node[above right] at (1.5,3) {$\y_{(3)}^{\quant}$};
\end{tikzpicture}}
    \caption{Quantitative poll $P_{(k)}^{\quant}$.}
    \label{subfig:MADS1}
  \end{subfigure}%
  \hfill
  \begin{subfigure}[t]{0.49\textwidth}
    \centering
    \captionsetup{justification=centering}
    \scalebox{1.05}{\begin{tikzpicture}
        \draw[step=0.125cm,gray, very thin] (0,0) grid (4,4);
        \draw[thick] (1.5,1.5) rectangle (2.5,2.5);

        \draw[dashed, -{Latex[length=3mm, width=2mm]}, thick] (2,2) -- (1.5,2);   
        \fill[black] (1.5,2) circle (2pt);
        \draw[dashed, -{Latex[length=3mm, width=2mm]}, thick] (2,2) -- (2,1.5);   
        \fill[black] (2,1.5) circle (2pt);
        \draw[dashed, -{Latex[length=3mm, width=2mm]}, thick] (2,2) -- (2,2.5);   
        \fill[black] (2,2.5) circle (2pt);
        \draw[dashed, -{Latex[length=3mm, width=2mm]}, thick] (2,2) -- (2.5,2);   
        \fill[black] (2.5,2) circle (2pt);

        \fill[black] (2,2) circle (2pt) node[above left, xshift=1pt] {$\x_{*}$};

        \node[left] at (1.5,2) {$\y_{(4)}^{\quant}$};
        \node[below] at (2,1.5) {$\y_{(5)}^{\quant}$};
        \node[above right] at (1.5,2.5) {$\y_{(6)}^{\quant}$};
        \node[above right] at (2.5,1.75) {$\y_{(7)}^{\quant}$};

\end{tikzpicture}}
    \caption{Quantitative poll $P_{(k+1)}^{\quant}$.}
    \label{subfig:MADS2}
  \end{subfigure}
\end{minipage}%
\hfill
%
\begin{minipage}[b]{0.32\textwidth}
  \centering
  \begin{subfigure}[t]{\textwidth}
    \centering
    \captionsetup{justification=centering}
    \scalebox{0.8}{\begin{tikzpicture}

\definecolor{myblue}{RGB}{70,130,255}
\definecolor{mypurple}{RGB}{186,85,211}
\definecolor{myred}{RGB}{240,60,60}
\definecolor{myorangegreen}{RGB}{210,180,80}
\definecolor{mygreen}{RGB}{80,180,80}

\def\xA{-0.75} \def\xB{1}
\def\xC{-0.5}  \def\xD{1.25}
\def\yBottom{-0.5}
\def\yTop{1.0}

\def\QuarterPlane{(\xA,\yBottom) -- (\xB,\yBottom) -- (\xD,\yTop) -- (\xC,\yTop) -- cycle;}

\def\XDiv{8}
\def\YDiv{8}

\newcommand{\DrawGridQuarter}{%
    \foreach \i in {0,...,\XDiv} {
        \pgfmathsetmacro{\t}{\i/\XDiv}
        \pgfmathsetmacro{\xb}{\xA + (\xB - \xA)*\t}
        \pgfmathsetmacro{\xt}{\xC + (\xD - \xC)*\t}
        \draw[black!40] (\xb,\yBottom) -- (\xt,\yTop);
    }
    \foreach \j in {0,...,\YDiv} {
        \pgfmathsetmacro{\y}{\yBottom + (\yTop - \yBottom)*\j/\YDiv}
        \pgfmathsetmacro{\xl}{\xA + (\xC - \xA)*(\y - \yBottom)/(\yTop - \yBottom)}
        \pgfmathsetmacro{\xr}{\xB + (\xD - \xB)*(\y - \yBottom)/(\yTop - \yBottom)}
        \draw[black!40] (\xl,\y) -- (\xr,\y);
    }
}

\def\Gap{1.6}

\def\PlaneColor#1{%
    \ifcase#1 myorangegreen\or
    myred\or
    mypurple\or
    myblue\or
    mygreen\fi
}

\newcommand{\GetGridPoint}[2]{%
    \pgfmathsetmacro{\t}{#1/\XDiv}
    \pgfmathsetmacro{\xb}{\xA + (\xB - \xA)*\t}
    \pgfmathsetmacro{\xt}{\xC + (\xD - \xC)*\t}
    \pgfmathsetmacro{\y}{\yBottom + (\yTop - \yBottom)*#2/\YDiv}
    \pgfmathsetmacro{\x}{\xb + (\xt - \xb)*( (\y - \yBottom)/(\yTop - \yBottom) )}
    (\x,\y)
}

    \foreach \i in {0,...,3} {
        \begin{scope}[shift={(0, \i*\Gap)}, xscale=1.75, yscale=1]
            \fill[\PlaneColor{\i}, opacity=0.4] \QuarterPlane;
            \begin{scope}
                \clip \QuarterPlane;
                \DrawGridQuarter
            \end{scope}
        \end{scope}
    }

    \draw[decorate, decoration={brace, amplitude=8pt}, thick]
    (-1.2, 1*\Gap-0.7) -- (-1.2, 3*\Gap+1.1)
    node[midway, left=6pt, align=center] {\large $\mathcal{N}(\x_{(k)}; 3)$};
    
    \begin{scope}[shift={(0, 3.4-0.3)}, xscale=1.5, yscale=1] 

      \def\ybot{-0.03}
      \def\ytop{0.3571+0.375}
      \def\xleft{-0.29}
      \def\xright{0.72} 
    
      \def\xslant{0.15}

      \draw[black, thick]
        (\xleft, \ybot) --
        (\xright, \ybot) --
        ({\xright + \xslant}, \ytop) --
        ({\xleft + \xslant}, \ytop) -- cycle;

      \fill[black] (0.29, 0.351) circle (2pt)
        node[below right, xshift=-3pt, yshift=1pt] {\large $\x_{(k)}$};
    
      \def\xk{0.29}
      \def\yk{0.351}
      \def\colw{0.44}
      \def\rowh{0.25}
      \def\xslant{0.0417}

      \fill[black] ({\xk + \colw}, \yk) circle (2pt);
      \draw[dashed,-{Latex[length=1.25mm, width=2mm]},thick]
        (\xk, \yk) -- ({\xk + \colw-0.05}, \yk);

      \fill[black] ({\xk - \colw + \xslant - 0.05}, {\yk + \rowh + 0.1}) circle (2pt);
      \draw[dashed,-{Latex[length=1.25mm, width=2mm]},thick]
        (\xk, \yk) -- ({\xk - \colw + \xslant - 0.01}, {\yk + \rowh + 0.075});

      \fill[black] ({\xk - \colw - \xslant - 0.1}, {\yk - \rowh - 0.1}) circle (2pt);
      \draw[dashed,-{Latex[length=1.25mm, width=2mm]},thick]
        (\xk, \yk) -- ({\xk - \colw - \xslant - 0.05}, {\yk - \rowh - 0.075});

      \fill[black] (0.29, -1.25) circle (2pt);
      \draw[dashed,-{Latex[length=1.25mm, width=2mm]},thick]
        (0.29, 0.351) -- (0.29, -1.2);

      \fill[black] (0.29, 1.95) circle (2pt);
      \draw[dashed,-{Latex[length=1.25mm, width=2mm]},thick]
        (0.29, 0.351) -- (0.29, 1.9);

    \end{scope}

\end{tikzpicture}}
    \caption{\catmads (unconstrained) $P_{(k)}^{\cat}$.}
    \label{subfig:CatMADS}
  \end{subfigure}
\end{minipage}

\caption{\blr 
(\subref{subfig:MADS1}) and (\subref{subfig:MADS2}) illustrate two consecutive
\mads quantitative polls at iterations $k$ and $k+1$, respectively, with
$\x_{(k)}=\x_{(k+1)}$.
(\subref{subfig:CatMADS}) illustrates both the quantitative and categorical polls
in \catmads.
}
\label{fig:MADS}
\end{figure}

Arrows emerging from an incumbent solution represent a positive basis.
Quantitative components produced are on the \textit{mesh} in gray and within the black square, called the \textit{frame}.
The frame delimits the region in which a quantitative poll can produce quantitative components.

In the example in \cref{subfig:CatMADS}, $x_{(k)}^{\cat}=\text{\sc red}$ and the number of components in the neighborhood is three. 
The neighborhood contains the $\text{{\sc red}}$ (incumbent), $\text{\sc yellow}$, and $\text{\sc purple}$ categories.
Each colored plane represents the quantitative domain associated with one categorical component.
Vertical arrows correspond to the categorical poll, while the quantitative poll is performed within the  $\text{{\sc red}}$  plane.

If both the search and poll steps fail to improve the objective, the optional step, called the \textit{extended poll}, may be launched at Step 3. 
This step revisits points in the categorical poll that almost improved the 
incumbent value $f(\x_{(k)})$ by calibrating their quantitative variables with additional evaluations.
{\blr
The extended poll is not required for the basic convergence guarantees of \catmads.
It provides stronger local optimality guarantees by accounting for interactions between categorical and quantitative variables, but requires additional evaluations.
A parameter $\xi$ controls the relative deterioration in the objective value that is accepted for triggering an extended poll.
Larger values of $\xi$ increase the number of extended polls and the evaluation cost, but strengthen the local optimality guarantees.
This creates a trade-off between evaluation cost and local optimality: see~\cite{catmads} for details.
}

{\blr 
The default value in the original \catmads paper is $\xi=0.05$~\cite{catmads}.
This value allows an objective deterioration of up to $5\%$ to trigger an extended poll.
For $\xi=0.05$, \catmads targets a
$[\cat-\continuous,\integer]$ local minimum $\x_{\star}$.
This type of local optimality is characterized with respect to a given categorical
neighborhood $\mathcal{N}(\x_{\star};m)$.
\cref{fig:catmads_local_minimum} provides an example illustrating such a local minimum.
The categorical variable has three components represented by the colored planes,
the integer variable takes values in $\{1,2,3\}$, and each plane contains a
two-dimensional continuous domain.
The local minimum $\x_{\star}$ is represented by $\star$.
The gray region around $\y$ illustrates the joint interaction between the
categorical and continuous variables.
After moving to a neighboring categorical component, local variations of the
continuous variables are also considered.
In contrast, the integer variable is treated separately, as illustrated by the
horizontal dotted lines connecting neighboring integer values.
This explains the notation $[\cat-\continuous,\integer]$, where the hyphen
denotes the joint interaction between categorical and continuous variables,
while the comma separates the integer variables.

Under the assumptions of the \catmads convergence analysis, this notion of
local optimality is established.
In particular, the objective function is assumed to be locally Lipschitz
continuous with respect to the continuous variables around $\x_{\star}$.
The distances defining the categorical neighborhoods must also eventually be
fixed, ensuring that the categorical neighborhoods no longer change during the
later iterations.
Finally, the objective function is assumed to be lower semicontinuous at
$\x_{\star}$ and upper semicontinuous at its categorical neighbors, limiting
abrupt changes in the objective when categorical components vary.
For more details, see Theorems~3 and~4 of the original \catmads paper~\cite{catmads}.

}

\begin{figure}[t!]
\centering

\begin{minipage}[b]{0.25\textwidth}
    \centering
    \scalebox{0.8}{\begin{tikzpicture}

\clip (-1.55, -0.1) rectangle (3.1, 7);

\newcommand{\FACE}[5]{ 
  \node[] (C) at (#1, #2, #3) {}; 
  \fill[#5,opacity=0.4] 
    ($(C) + (#4*0.20,#4,#4)$) -- 
    ($(C) + (0,#4,#4)$) -- 
    ($(C) + (0,#4,0)$) -- 
    ($(C) + (#4*0.20,#4,0)$) -- 
    cycle;
  \draw[-] ($(C) + (#4*0.20,#4,#4)$) -- ($(C) + (0,#4,#4)$);
  \draw[-] ($(C) + (#4*0.20,#4,#4)$) -- ($(C) + (#4*0.20,#4,0)$);
  \draw[-] ($(C) + (0,#4,0)$) -- ($(C) + (0,#4,#4)$);
  \draw[-] ($(C) + (0,#4,0)$) -- ($(C) + (#4*0.20,#4,0)$);
}

\newcommand{\FaceSize}{4} 

\definecolor{blueface}{RGB}{51,171,255}
\definecolor{redface}{RGB}{236,87,87}
\definecolor{greenface}{RGB}{112,219,112}

\def\xleft{0}
\def\xmid{1.15}
\def\xright{2.3}

\FACE{\xleft}{2}{0}{\FaceSize}{greenface}
\FACE{\xleft}{0}{0}{\FaceSize}{redface}
\FACE{\xleft}{-2}{0}{\FaceSize}{blueface}
\fill[black] ($(\xleft,0,0) + (\FaceSize*0.10,\FaceSize, \FaceSize/2)$) circle (2pt);

\FACE{\xmid}{2}{0}{\FaceSize}{greenface}
\fill[black] ($(\xmid,2,0) + (\FaceSize*0.10,\FaceSize, \FaceSize/2)$) circle (2pt);
\filldraw[gray, thick, rotate around={-45:(\xmid+\FaceSize*0.10, 2+\FaceSize, \FaceSize/2)}] 
  (\xmid+\FaceSize*0.10, 2+\FaceSize, \FaceSize/2) ellipse (6.5pt and 11pt);
\FACE{\xmid}{0}{0}{\FaceSize}{redface}
\filldraw[gray, thick, rotate around={-45:(\xmid+\FaceSize*0.10, 0+\FaceSize, \FaceSize/2)}] 
  (\xmid+\FaceSize*0.10, 0+\FaceSize, \FaceSize/2) ellipse (6.5pt and 11pt);
\node at ($(\xmid,0,0) + (\FaceSize*0.10,\FaceSize, \FaceSize/2)$) {\LARGE $\star$};
\FACE{\xmid}{-2}{0}{\FaceSize}{blueface}

\FACE{\xright}{2}{0}{\FaceSize}{greenface}
\FACE{\xright}{0}{0}{\FaceSize}{redface}
\fill[black] ($(\xright,0,0) + (\FaceSize*0.10,\FaceSize, \FaceSize/2)$) circle (2pt);
\FACE{\xright}{-2}{0}{\FaceSize}{blueface}

\draw[dashed,-{Latex[length=1.25mm, width=2mm]},thick] 
  ($(\xmid,0,0) + (\FaceSize*0.10,\FaceSize, \FaceSize/2)$) -- 
  ($(\xmid,2,0) + (\FaceSize*0.10,\FaceSize, \FaceSize/2)$);

\draw[dashed,-{Latex[length=1.25mm, width=2mm]},thick] 
  ($(\xmid,0,0) + (\FaceSize*0.10,\FaceSize, \FaceSize/2)$) -- 
  ($(\xleft+0.1,0,0) + (\FaceSize*0.10,\FaceSize, \FaceSize/2)$);

\draw[dashed,-{Latex[length=1.25mm, width=2mm]},thick] 
  ($(\xmid,0,0) + (\FaceSize*0.10,\FaceSize, \FaceSize/2)$) -- 
  ($(\xright,0,0) + (\FaceSize*0.10-0.1,\FaceSize, \FaceSize/2)$);

\node[anchor=south, text=gray] at ($(\xmid+0.5,2.25,0) + (\FaceSize*0.10,\FaceSize, \FaceSize/2)$) {\Large $\bm{y}$};
\node[anchor=west] at ($(\xright,0,0) + (\FaceSize*0.10,\FaceSize, \FaceSize/2)$) {\Large $\z$};

\node[] at (0.5, 6.5, 0) {\Large $1$};
\node[] at (1.6, 6.5, 0) {\Large $2$};
\node[] at (2.7, 6.5, 0) {\Large $3$};

\end{tikzpicture}}
    \captionsetup{justification=centering}
    \caption{$[\cat-\continuous,\integer]$ local minimum in \catmads \hspace{0.25cm} (reproduced from~\cite{catmads}).}
    \label{fig:catmads_local_minimum}
\end{minipage}
\hfill
\begin{minipage}[b]{0.675\textwidth}
    \centering

    \nextfloat

    \begin{subfigure}[b]{0.47\textwidth}
      \centering
      \scalebox{0.9}{\newcommand{\curveGreen}{
      (-2.2,-1)
      .. (0.2,-1.5)
      .. (0.8,-1.7) 
      .. (2,2)
      .. (0,1)
}

\newcommand{\curveBlue}{(-1,0) .. (0.5,3) .. (2,2) .. (0,1)}

\newcommand{\curveRed}{
  (-2,-0.8)
  .. (0.2,-1.5)
  .. (1.0,-0.8)
  .. (1.5,2)
  .. (0.5,1.2)
  .. (0,1)
}

\begin{tikzpicture}[use Hobby shortcut]
\definecolor{myblue}{RGB}{70,130,255}
\definecolor{mypurple}{RGB}{186,85,211}
\definecolor{myred}{RGB}{240,60,60}
\definecolor{myorangegreen}{RGB}{210,180,80}
\definecolor{mygreen}{RGB}{80,180,80}

\def\xA{-0.75} \def\xB{1}
\def\xC{-0.5}  \def\xD{1.25}
\def\yBottom{-0.5}
\def\yTop{1.0}

\def\QuarterPlane{(\xA,\yBottom) -- (\xB,\yBottom) -- (\xD,\yTop) -- (\xC,\yTop) -- cycle;}

\def\XDiv{8}
\def\YDiv{8}

\newcommand{\DrawGridQuarter}{%
    \foreach \i in {0,...,\XDiv} {
        \pgfmathsetmacro{\t}{\i/\XDiv}
        \pgfmathsetmacro{\xb}{\xA + (\xB - \xA)*\t}
        \pgfmathsetmacro{\xt}{\xC + (\xD - \xC)*\t}
        \draw[black!40] (\xb,\yBottom) -- (\xt,\yTop);
    }
    \foreach \j in {0,...,\YDiv} {
        \pgfmathsetmacro{\y}{\yBottom + (\yTop - \yBottom)*\j/\YDiv}
        \pgfmathsetmacro{\xl}{\xA + (\xC - \xA)*(\y - \yBottom)/(\yTop - \yBottom)}
        \pgfmathsetmacro{\xr}{\xB + (\xD - \xB)*(\y - \yBottom)/(\yTop - \yBottom)}
        \draw[black!40] (\xl,\y) -- (\xr,\y);
    }
}

\def\Gap{1.6}

\def\PlaneColor#1{%
    \ifcase#1 myblue\or
    mypurple\or
    myred\or
    myorangegreen\or
    mygreen\fi
}

\newcommand{\GetGridPoint}[2]{%
    \pgfmathsetmacro{\t}{#1/\XDiv}
    \pgfmathsetmacro{\xb}{\xA + (\xB - \xA)*\t}
    \pgfmathsetmacro{\xt}{\xC + (\xD - \xC)*\t}
    \pgfmathsetmacro{\y}{\yBottom + (\yTop - \yBottom)*#2/\YDiv}
    \pgfmathsetmacro{\x}{\xb + (\xt - \xb)*( (\y - \yBottom)/(\yTop - \yBottom) )}
    (\x,\y)
}

\foreach \i in {1,...,4} {
    \begin{scope}[shift={(0, \i*\Gap)}, xscale=1.75, yscale=1]
        \fill[gray, opacity=0.25] \QuarterPlane;
        \begin{scope}
            \clip \QuarterPlane;
            \DrawGridQuarter
        \end{scope}
    \end{scope}
}


\begin{scope}[shift={(0, 4*\Gap)}, xscale=1.75, yscale=1]
  \begin{scope}[scale=0.2, rotate=30, shift={(0.25, 0.75)}]
    \draw[closed, scale=1, fill=myblue, draw=black!40, opacity=0.6] \curveGreen;
  \end{scope}
\end{scope}

\begin{scope}[shift={(0, 1*\Gap)}, xscale=1.75, yscale=1]
  \begin{scope}[scale=0.275, rotate=100, shift={(1, -3.75)}]
    \draw[use Hobby shortcut, closed, fill=myorangegreen, draw=black!40, opacity=0.6] \curveBlue;
  \end{scope}
\end{scope}

\begin{scope}[shift={(0, 3*\Gap)}, xscale=1.75, yscale=1]
  \begin{scope}[scale=0.2, rotate=-30, shift={(-0.75, -0.5)}]
    \draw[use Hobby shortcut, closed, fill=mypurple, draw=black!40, opacity=0.6] \curveRed;
  \end{scope}
\end{scope}

\begin{scope}[shift={(0, 2*\Gap)}, xscale=1.75, yscale=1]
  \begin{scope}[scale=0.4, rotate=30, shift={(0.75, 0.25)}]
    \fill[myred, opacity=0.6, draw=black!40]
      (0, 0) ellipse [x radius=2, y radius=0.8];
  \end{scope}
\end{scope}


\begin{scope}[shift={(0-0.375, 3.4-0.3)}, xscale=1.5, yscale=1]

  \def\ybot{-0.03}
  \def\ytop{0.3571+0.375}
  \def\xleft{-0.29}
  \def\xright{0.72}
  \def\xslant{0.15}

  \draw[black]
      (\xleft, \ybot) --
      (\xright, \ybot) --
      ({\xright + \xslant}, \ytop) --
      ({\xleft + \xslant}, \ytop) -- cycle;

    \fill[black] (0.29, 0.351) circle (2pt)
    node[above, xshift=4pt, yshift=0pt] {$\x_{\mathrm{FEA}}$};

    %

  \def\xk{0.29}
  \def\yk{0.351}
  \def\colw{0.44}
  \def\rowh{0.25}
  \def\xslant{0.0417}

  \fill[black] (0.29, -1.25) circle (2pt);
  \draw[dashed,-{Latex[length=1.25mm, width=2mm]},thick]
    (0.29, 0.351) -- (0.29, -1.2);

\end{scope}

\begin{scope}[shift={(0+0.77, 3.4-0.3+1.6)}, xscale=1.5, yscale=1]

  \def\ybot{-0.03}
  \def\ytop{0.3571+0.375}
  \def\xleft{-0.29}
  \def\xright{0.72}
  \def\xslant{0.15}

  \draw[black]
      (\xleft, \ybot) --
      (\xright, \ybot) --
      ({\xright + \xslant}, \ytop) --
      ({\xleft + \xslant}, \ytop) -- cycle;

  \fill[black] (0.29, 0.351) circle (2pt)
    node[above, xshift=-1pt, yshift=0pt] {$\x_{\mathrm{INF}}$};

  \def\xk{0.29}
  \def\yk{0.351}
  \def\colw{0.44}
  \def\rowh{0.25}
  \def\xslant{0.0417}

  \fill[black] (0.29, 1.95) circle (2pt);
  \draw[dashed,-{Latex[length=1.25mm, width=2mm]},thick]
    (0.29, 0.351) -- (0.29, 1.9);

    %

\end{scope}

\end{tikzpicture}}
      \caption{Basic categorical poll in~\cite{catmads}.}
      \label{subfig:CatMADS_constrained1}
    \end{subfigure}
    \hfill
    \begin{subfigure}[b]{0.52\textwidth}
      \centering
      \scalebox{0.9}{\newcommand{\curveGreen}{
      (-2.2,-1)
      .. (0.2,-1.5)
      .. (0.8,-1.7) 
      .. (2,2)
      .. (0,1)
}

\newcommand{\curveBlue}{(-1,0) .. (0.5,3) .. (2,2) .. (0,1)}

\newcommand{\curveRed}{
  (-2,-0.8)
  .. (0.2,-1.5)
  .. (1.0,-0.8)
  .. (1.5,2)
  .. (0.5,1.2)
  .. (0,1)
}

\begin{tikzpicture}[use Hobby shortcut]
\definecolor{myblue}{RGB}{70,130,255}
\definecolor{mypurple}{RGB}{186,85,211}
\definecolor{myred}{RGB}{240,60,60}
\definecolor{myorangegreen}{RGB}{210,180,80}
\definecolor{mygreen}{RGB}{80,180,80}

\def\xA{-0.75} \def\xB{1}
\def\xC{-0.5}  \def\xD{1.25}
\def\yBottom{-0.5}
\def\yTop{1.0}

\def\QuarterPlane{(\xA,\yBottom) -- (\xB,\yBottom) -- (\xD,\yTop) -- (\xC,\yTop) -- cycle;}

\def\XDiv{8}
\def\YDiv{8}

\newcommand{\DrawGridQuarter}{%
    \foreach \i in {0,...,\XDiv} {
        \pgfmathsetmacro{\t}{\i/\XDiv}
        \pgfmathsetmacro{\xb}{\xA + (\xB - \xA)*\t}
        \pgfmathsetmacro{\xt}{\xC + (\xD - \xC)*\t}
        \draw[black!40] (\xb,\yBottom) -- (\xt,\yTop);
    }
    \foreach \j in {0,...,\YDiv} {
        \pgfmathsetmacro{\y}{\yBottom + (\yTop - \yBottom)*\j/\YDiv}
        \pgfmathsetmacro{\xl}{\xA + (\xC - \xA)*(\y - \yBottom)/(\yTop - \yBottom)}
        \pgfmathsetmacro{\xr}{\xB + (\xD - \xB)*(\y - \yBottom)/(\yTop - \yBottom)}
        \draw[black!40] (\xl,\y) -- (\xr,\y);
    }
}

\def\Gap{1.6}

\def\PlaneColor#1{%
    \ifcase#1 myblue\or
    mypurple\or
    myred\or
    myorangegreen\or
    mygreen\fi
}

\newcommand{\GetGridPoint}[2]{%
    \pgfmathsetmacro{\t}{#1/\XDiv}
    \pgfmathsetmacro{\xb}{\xA + (\xB - \xA)*\t}
    \pgfmathsetmacro{\xt}{\xC + (\xD - \xC)*\t}
    \pgfmathsetmacro{\y}{\yBottom + (\yTop - \yBottom)*#2/\YDiv}
    \pgfmathsetmacro{\x}{\xb + (\xt - \xb)*( (\y - \yBottom)/(\yTop - \yBottom) )}
    (\x,\y)
}

\foreach \i in {1,...,4} {
    \begin{scope}[shift={(0, \i*\Gap)}, xscale=1.75, yscale=1]
        \fill[gray, opacity=0.25] \QuarterPlane;
        \begin{scope}
            \clip \QuarterPlane;
            \DrawGridQuarter
        \end{scope}
    \end{scope}
}


\begin{scope}[shift={(0, 4*\Gap)}, xscale=1.75, yscale=1]
  \begin{scope}[scale=0.2, rotate=30, shift={(0.25, 0.75)}]
    \draw[closed, scale=1, fill=myblue, draw=black!40, opacity=0.6] \curveGreen;
  \end{scope}
\end{scope}

\begin{scope}[shift={(0, 1*\Gap)}, xscale=1.75, yscale=1]
  \begin{scope}[scale=0.275, rotate=100, shift={(1, -3.75)}]
    \draw[use Hobby shortcut, closed, fill=myorangegreen, draw=black!40, opacity=0.6] \curveBlue;
  \end{scope}
\end{scope}

\begin{scope}[shift={(0, 3*\Gap)}, xscale=1.75, yscale=1]
  \begin{scope}[scale=0.2, rotate=-30, shift={(-0.75, -0.5)}]
    \draw[use Hobby shortcut, closed, fill=mypurple, draw=black!40, opacity=0.6] \curveRed;
  \end{scope}
\end{scope}

\begin{scope}[shift={(0, 2*\Gap)}, xscale=1.75, yscale=1]
  \begin{scope}[scale=0.4, rotate=30, shift={(0.75, 0.25)}]
    \fill[myred, opacity=0.6, draw=black!40]
      (0, 0) ellipse [x radius=2, y radius=0.8];
  \end{scope}
\end{scope}


\begin{scope}[shift={(0-0.375, 3.4-0.3)}, xscale=1.5, yscale=1]

  \def\ybot{-0.03}
  \def\ytop{0.3571+0.375}
  \def\xleft{-0.29}
  \def\xright{0.72}
  \def\xslant{0.15}

  \draw[black]
      (\xleft, \ybot) --
      (\xright, \ybot) --
      ({\xright + \xslant}, \ytop) --
      ({\xleft + \xslant}, \ytop) -- cycle;

    \fill[black] (0.29, 0.351) circle (2pt)
    node[above, xshift=4pt, yshift=0pt] {$\x_{\mathrm{FEA}}$};

    \fill[black] (0.29, 1.95+1.6) circle (2pt); 
    \draw[dashed,-{Latex[length=1.25mm, width=2mm]},thick] (0.29, 0.351) -- (0.29, 1.9+1.6);

  \def\xk{0.29}
  \def\yk{0.351}
  \def\colw{0.44}
  \def\rowh{0.25}
  \def\xslant{0.0417}


\end{scope}

\begin{scope}[shift={(0+0.77, 3.4-0.3+1.6)}, xscale=1.5, yscale=1]

  \def\ybot{-0.03}
  \def\ytop{0.3571+0.375}
  \def\xleft{-0.29}
  \def\xright{0.72}
  \def\xslant{0.15}

  \draw[black]
      (\xleft, \ybot) --
      (\xright, \ybot) --
      ({\xright + \xslant}, \ytop) --
      ({\xleft + \xslant}, \ytop) -- cycle;

  \fill[black] (0.29, 0.351) circle (2pt)
    node[above, xshift=-1pt, yshift=0pt] {$\x_{\mathrm{INF}}$};

  \def\xk{0.29}
  \def\yk{0.351}
  \def\colw{0.44}
  \def\rowh{0.25}
  \def\xslant{0.0417}


    \fill[black] (0.29, -1.25) circle (2pt); 
    \draw[dashed,-{Latex[length=1.25mm, width=2mm]},thick] (0.29, 0.351) -- (0.29, -1.2);

\end{scope}

\end{tikzpicture}}
      \caption{Poll with surrogate-based neighborhoods.}
      \label{subfig:CatMADS_surrogate}
    \end{subfigure}


    \caption{The \catmads framework in the presence of constraints.}
    \label{fig:CatMADS_pair}
\end{minipage}

\end{figure}

{\blr In} an unsuccessful iteration, the mesh and frame are further discretized.
\catmads terminates when the mesh and frame have reached their minimum allowable discretization.


In \catmads, the constraints are handled with the progressive barrier strategy (PB)~\cite{AuDe09a}.
{\blr
The PB considers infeasible trial points and progressively reduces their constraint violation.
For inequality constraints $g_j(\x)\leq 0$, $j\in J$, the constraint violation is measured by the aggregation function
\begin{equation*}
    h(\x)\coloneq
    \sum_{j\in J}\left(\max\{0,g_j(\x)\}\right)^2.
\end{equation*}
A point is feasible if and only if $h(\x)=0$.
}
The PB works with two incumbent solutions, each with their own independent poll, as described above.
The \textit{feasible poll} is performed around the \textit{feasible incumbent solution} $\smash{\x_{\feasible}}$, and this solution is {\blr replaced when a feasible point} strictly improves the objective function.

The \textit{infeasible poll} is performed at the \textit{infeasible incumbent solution} $\smash{\x_{\infeasible}}$, and this solution is forced to become more feasible through the iterations.
{\blr For more details on the update with the PB, see~\cite[Chapter~12]{AuHa2026}~or~\cite{AuDe09a}.}
%

\subsection{Improvements of \catmads with surrogate-based neighborhoods}
\label{sec:CatMADS_improved}

In~\cite{catmads}, categorical neighborhoods are constructed using information
from the objective function only.
Moreover, the implementation in~\cite{catmads} relies on a simple mixed-variable
interpolation.
The surrogate-based neighborhoods introduced here do not affect the theoretical
convergence results, as the proposed method remains an instance of the
\catmads\ framework.
However, the resulting categorical neighborhoods incorporate information from
both the objective and the constraint functions and rely on an objective
pseudo-distance induced by kernel-based similarities.
As a result, the categorical polls should be guided more effectively toward feasible regions than in the original implementation.
Again, this has no impact on the theoretical results, but it is expected to improve empirical performance and accelerate convergence.
The improved categorical poll with the surrogate-based neighborhoods is illustrated in \cref{fig:CatMADS_pair}.

In \cref{subfig:CatMADS_constrained1}, the infeasible regions are in gray and the feasible regions are colored.
The feasible categorical poll starts from the red category and selects the yellow category, which is the most similar with respect to the objective function.
The generated trial point is infeasible.

In contrast, in \cref{subfig:CatMADS_surrogate}, the generated trial point lies in the blue category.
%
%
%

A similar behavior is observed for the infeasible categorical poll.
In the basic implementation, the generated trial point comes from the category that is closest with respect to the objective function, which again leads to an infeasible trial point in \cref{subfig:CatMADS_constrained1}.
In comparison, the surrogate-based implementation selects a category that is less similar in terms of the objective but considers the similarities of both the objective and the constraint functions.
This results in a feasible trial point in the red category in \cref{subfig:CatMADS_surrogate}.

{\blr 
The \catmads convergence theory remains valid if the pseudo-distance functions are fixed after a finite number of iterations.
This condition requires the Gaussian process models to stop being updated after some iteration.
The categorical poll otherwise follows the standard \catmads construction, with the
quantitative component fixed and the categorical components selected from the
surrogate-based neighborhood.
The local optimality conditions are then defined with respect to the
surrogate-based neighborhoods.
}

\section{Computational experiments}
\label{sec:computational_experiments}

This section benchmarks an instance of \catmads using the novel neighborhoods described in \cref{sec:direct_search}.
This is a proof of concept demonstrating the utility of surrogate-based neighborhoods in mixed-variable blackbox optimization.
Since the surrogate models are GPs with one-hot encoding of the categorical variables, the new resulting method is referred to as \catmadsgp.

{\blr
The numerical experiments are organized as follows.
First, the implementation of \catmadsgp is presented in \cref{sec:implementation_details}.
The data profiles used for benchmarking are then introduced in \cref{sec:data_profiles}.
The first comparison, presented in \cref{sec:data_profiles_catmads}, introduces two additional \catmads variants and compares them with \catmadsgp to isolate the contributions of the proposed neighborhood construction.
The second comparison, presented in \cref{sec:data_profiles_solvers}, compares \catmadsgp with external mixed-variable solvers.
Finally, the mechanical-part problem is used as an application in \cref{sec:A}.
}

\subsection{Implementation of the \catmads variants}
\label{sec:implementation_details}

{\blr
%
%
%
%
%
%
%
%

The new surrogate-based neighborhoods are benchmarked within the \catmads framework.
\catmads is well suited for this comparison because its categorical poll is defined from categorical neighborhoods.
Three variants are considered.
The first follows the original \catmads neighborhood construction~\cite{catmads} and is denoted \catmadsf.
In this variant, the categorical neighborhoods are constructed from a mixed-variable interpolation model of the objective function.
The interpolation model is selected by cross-validation and uses only objective information.
%
%
For the second variant, the categorical neighborhoods remain based only on the objective function.
However, the mixed-variable interpolation model is replaced by a mixed-variable Gaussian process.
This variant is denoted \catmadsgpf.
The third variant is \catmadsgp: it uses surrogate-based neighborhoods constructed from Gaussian process models of both the objective and constraint functions.

All three variants are implemented within \nomadbbo, a Python package built on NOMAD~4.6 and currently available as a beta release at \url{https://test.pypi.org/project/NomadBBO/}.
They follow the same optimization procedure and differ only in the construction of the categorical neighborhoods.
These neighborhoods determine the points generated by the categorical poll.
The common optimization procedure is summarized in \cref{algo:catmadsgp}.

\begin{algorithm}[htb!]
\SetAlgoLined

Given $f:\mathcal{X}\to\overline{\mathbb{R}}$, constraint functions
$g_j:\mathcal{X}\to\overline{\mathbb{R}}$, $j\in J$, and an evaluation budget.

\SetKwBlock{init}{\textbf{0. Initialization}}{}
\init{
Perform an initial DoE and initialize the \catmads framework\;

Construct surrogate model(s) from the initial DoE\;

}

\SetKwFor{While}{While}{do}{end}

\While{no stopping criterion is met}{

    \textbf{1. (Optional) Opportunistic search.}\;
    \hspace{0.33cm} Perform the standard \catmads search step\;

    \vspace{0.05cm}

    \textbf{2. Opportunistic poll.}\;
    %
    %
    \hspace{0.33cm} Construct the surrogate-based neighborhoods from current incumbents \;
    \hspace{0.33cm}  Perform the infeasible and feasible categorical polls from surrogate-based neighborhoods\;
    \hspace{0.33cm} Perform the infeasible and feasible quantitative {\blr polls}\;

    \vspace{0.05cm}

    \textbf{3. (Optional) Opportunistic extended poll.}\;
    \hspace{0.33cm} Perform the extended poll when applicable\;

    \vspace{0.05cm}

    \textbf{4. Update.}\;
    \hspace{0.33cm} Update the incumbent solutions, progressive barrier, and mesh parameters of \catmads\;
    \hspace{0.33cm} Update the set of evaluated points\;
}

\caption{\catmads framework for the three variants.}
\label{algo:catmadsgp}
\end{algorithm}

For \catmadsgp, the categorical neighborhoods are constructed from Gaussian process
models of the objective and constraint functions.
These models are implemented with the {\sf SMT Python library}, version~2.15.0~\cite{BaBuDiHwMaMoLaLeSa2023}, using one-hot encoding with continuous relaxation for the categorical variables and the default squared exponential kernel for the quantitative variables.
All other model parameters are left at their default values in {\sf SMT}~2.15.0.
Non-finite output values are excluded when fitting the models.
The models are constructed once at initialization using the initial DoE and remain fixed throughout the optimization.
This ensures that the pseudo-distance functions remain fixed, as required by the \catmads convergence assumptions.

The normalization bounds used in $\psi$ are determined from the minimum and maximum finite constraint values in the initial DoE and are kept fixed throughout the optimization.
If the normalization range is zero, the normalized values are set to $1$.
For the constraint pseudo-distance, the simplest parameter choices are used.
The relaxation parameter is set to $\lambda=0$ and the Euclidean norm is used, \textit{i.e.}, $p=2$.
%
%
}

{\blr
For a given problem, the number of neighbors is
$m=\max\left\{3,\sqrt{\left|\mathcal{X}^{\cat}\right|}\right\}$,
rounded to the nearest integer~\cite{catmads}.
Ties in the categorical ordering are resolved according to the order in
which the components are encountered when enumerating
$\mathcal{X}^{\cat}$.
}



{\blr

\subsection{Data profiles}
\label{sec:data_profiles}

The main benchmarking tools used in this study are data profiles.
Data profiles compare the performance of several solvers over a set of test problems.
The problems considered are taken from the \catsuite collection~\cite{catsuite}.
The collection contains $30$ unconstrained and $30$ constrained mixed-variable problems.

Each problem is instantiated with five different seeds, resulting in a total of 300 instances. 
An instance, denoted $p \in \mathcal{P}$, corresponds to a problem with a specific seed.
For each instance, all solvers use the same initial DoE generated with the corresponding seed.
For a given data profile, the set of considered instances is denoted $\mathcal{P}$.
The convergence test for an instance $p \in \mathcal{P}$ depends on an initial objective value $f_0 \in \mathbb{R}$.
%
For the unconstrained problems, $f_0$ is set to the smallest objective function value in the common DoE. In the constrained case, $f_0$ is set to the smallest of the first feasible objective function values found by the solvers~\cite{G-2025-36}.
This definition enables the comparison of solvers even when the DoE fails to produce a feasible solution.

Each problem is assigned a fixed evaluation budget proportional to \(n_p+1\), where \(n_p\) denotes the number of variables of the problem. 
The specific budget is given in the corresponding subsections.
%

%
The set of solvers considered in a comparison is denoted by \(\mathcal{S}\).
For example, in \cref{sec:data_profiles_catmads}, the three \catmads variants are compared, such that \(\mathcal{S}=\{ \catmadsgp, \catmadsgpf, \catmadsf \}\).
A solver $s \in \mathcal{S}$ $\tau$-solves an instance $p \in \mathcal{P}$ if it produces a feasible incumbent solution $\x_{\feasible} \in \Omega$ such that the reduction $f_0 - f(\x_{\feasible})$ is within $\tau$ of the largest reduction $f_0 - f_{\star}$ obtained by any solver on this instance~\cite{MoWi2009}.
More precisely, the $\tau$-convergence test is defined as
\begin{equation*}
f_0 - f(\x_{\feasible}) \ge (1-\tau)(f_0 - f_{\star}),
\end{equation*} 
where $\tau \in [0,1]$ is a given tolerance.

A data profile~\cite{MoWi2009} represents the fraction of instances $\tau$-solved by a solver $s \in \mathcal{S}$ as a function of the computational budget. 
It is defined as
\begin{equation*}
    \text{data}_s\left( \alpha \right) \coloneq \frac{1}{\left| \mathcal{P} \right|} \left| \left\{ p \in \mathcal{P} \, : \, \frac{k_{p, s}}{n_p + 1} \leq \alpha \right\} \right| \in [0,1],
\end{equation*}
where $k_{p,s} \geq 0$ is the number of evaluations required by solver $s$ to $\tau$-solve instance $p$, and $\alpha \geq 0$ specifies the evaluation budget in multiples of $(n_p+1)$ evaluations.
{\blr For a given value of $\alpha$, solver $s$ must solve instance $p$ within $\alpha(n_p+1)$ evaluations.}

\subsection{Comparison of the \catmads neighborhood strategies}
\label{sec:data_profiles_catmads}

In this section, the surrogate-based neighborhoods are benchmarked within the \catmads framework using the three variants introduced in \cref{sec:implementation_details}.
For this experiment, the optional search and extended poll from \cref{algo:catmadsgp} are disabled for all three variants, since the comparison focuses on the contribution of the categorical neighborhoods.
The variants therefore use the same optimization mechanisms except for the neighborhoods used to generate the categorical poll points.
The evaluation budget is set to $100 (n_p+1)$ for each problem.

The data profiles for the three variants are presented in \cref{fig:data_profiles_CatMADS}.
For the unconstrained problems, \catmadsgp and \catmadsgpf are equivalent because there are no constraints and the categorical neighborhoods use only the Gaussian process of the objective function.

\begin{figure}[htb!]
\centering
\begin{subfigure}{\textwidth}
\centering
\includegraphics[width=1\linewidth]{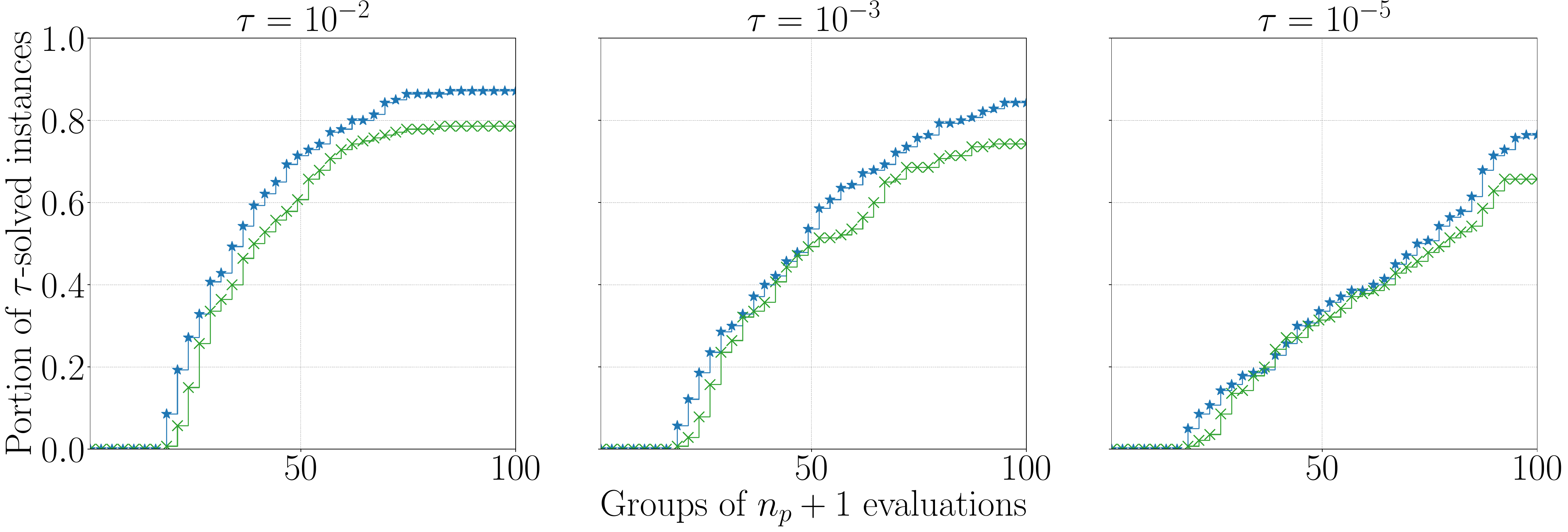}
\caption{Unconstrained test problems.}
\label{subfig:dataprofiles_unconstrained_CatMADS}
\end{subfigure}
\begin{subfigure}{\textwidth}
\centering
\includegraphics[width=1\linewidth]{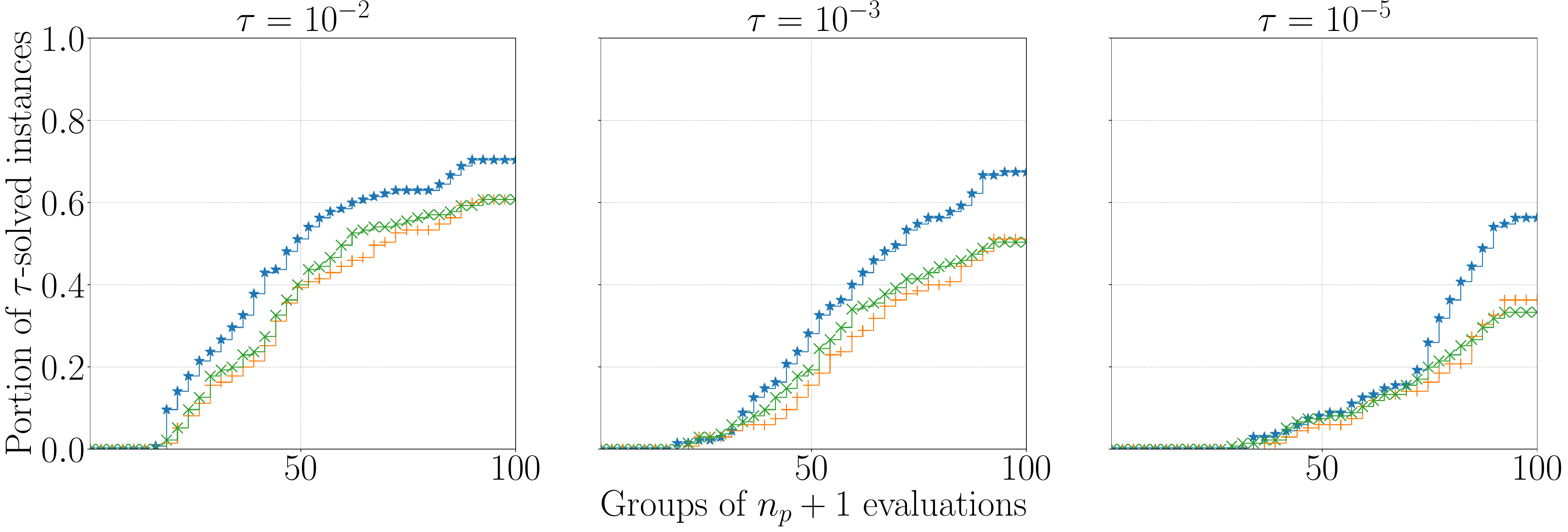}
\caption{Constrained test problems.}
\label{subfig:dataprofiles_constrained_CatMADS}
\end{subfigure}
\begin{subfigure}{\textwidth}
\centering
\includegraphics[width=0.55\linewidth]{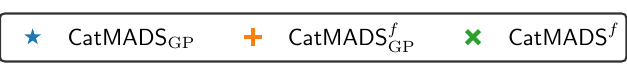}
\end{subfigure}
\caption{Comparison results between the \catmads variants. }
\label{fig:data_profiles_CatMADS}
\end{figure}

For the unconstrained data profile, \catmadsgp performs better than \catmadsf for all three tolerances.
This suggests an improvement from the (unconstrained) surrogate-based neighborhoods over the original \catmadsf neighborhoods on these problems.

For the constrained data profiles, \catmadsgp yields the best performance.
The difference is most pronounced for the strictest tolerance, $\tau=10^{-5}$.
For instance, at $100(n_p+1)$ evaluations, \catmadsgp $\tau$-solves approximately $57\%$ of the instances, compared with about $36\%$ for \catmadsgpf and $33\%$ for \catmadsf.
The similar performance of \catmadsgpf and \catmadsf suggests that the main improvement for the \catmads framework comes from the surrogate-based neighborhoods incorporating constraint information.

}

\clearpage \newpage

\subsection{Comparison with existing solvers}
\label{sec:external_solvers}

{\blr
The next experiments compare \catmadsgp with two external mixed-variable solvers used in blackbox optimization applications~\cite{JiKu2023,PrTr2024,SrArAh2025}.
Two experiments are considered.
The first in \cref{sec:data_profiles_solvers}  uses data profiles on the \catsuite test problems to compare the overall performance of the solvers.
The second in \cref{sec:A} considers the mechanical-part application over 50 independent runs.

The two external mixed-variable solvers considered are:
\begin{itemize}

    \item \href{https://pymoo.org/customization/mixed.html}{\sf Pymoo}~0.6.1.6, a Python library containing a collection of optimization algorithms.
    Specifically, the mixed-variable implementation of a genetic algorithm is used~\cite{BlankDeb2020};

    \item \href{https://hub.optuna.org/tags/cma-es/}{\sf Optuna}~5.0, a mixed-variable solver extending the covariance matrix adaptation evolution strategy (CMA-ES) to categorical variables~\cite{AkSaYaOhKo2020}.
    
\end{itemize}
{\sf Pymoo} and {\sf Optuna} are used with the default settings of their mixed-variable solvers.
The initial points and evaluation budget are specified so that all solvers use the same values in a given experiment.
}

\subsubsection{Data profiles}
\label{sec:data_profiles_solvers}

{\blr
This comparison focuses on the overall performance of the solvers rather than only on differences between \catmads variants.
A larger budget of \(250 (n_p+1)\) is used to allow global or more exploratory methods sufficient evaluations to identify good solutions.
The same budget and initial DoE are used for all solvers, with the DoE representing \(20\%\) of the budget.
}

\begin{figure}[htb!]
\centering
\begin{subfigure}{\textwidth}
\centering
\includegraphics[width=1\linewidth]{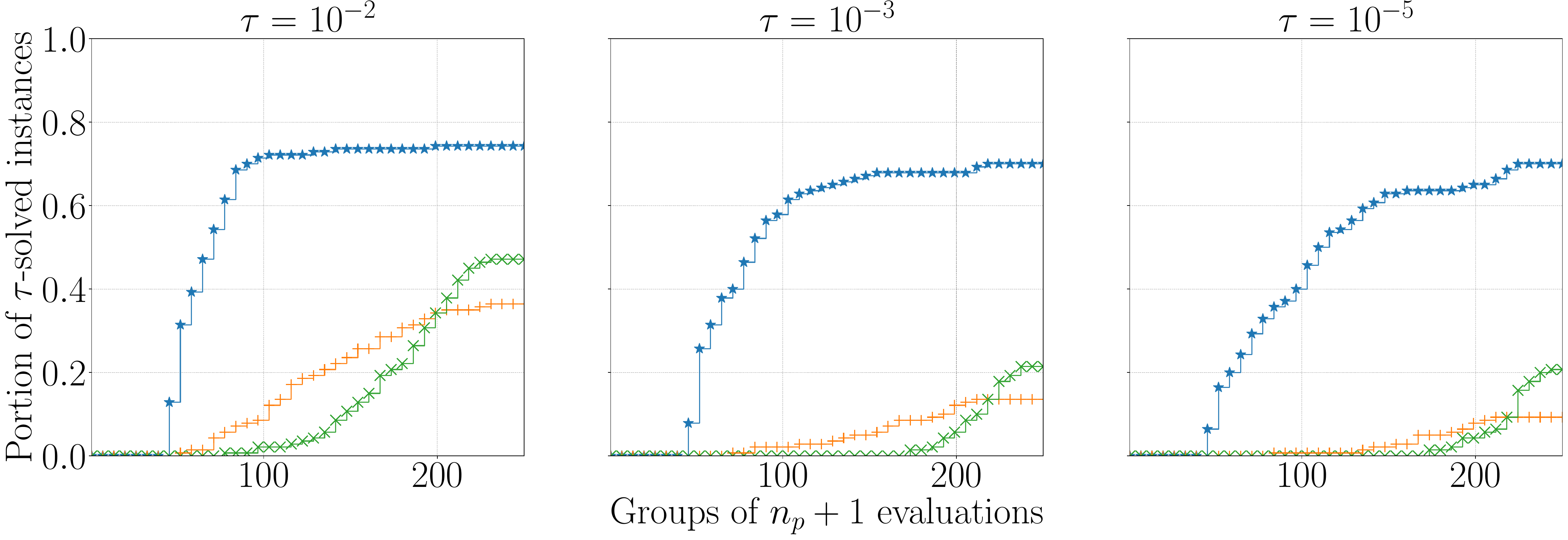}
\caption{Unconstrained test problems.}
\label{subfig:dataprofiles_unconstrained}
\end{subfigure}
\begin{subfigure}{\textwidth}
\centering
\includegraphics[width=1\linewidth]{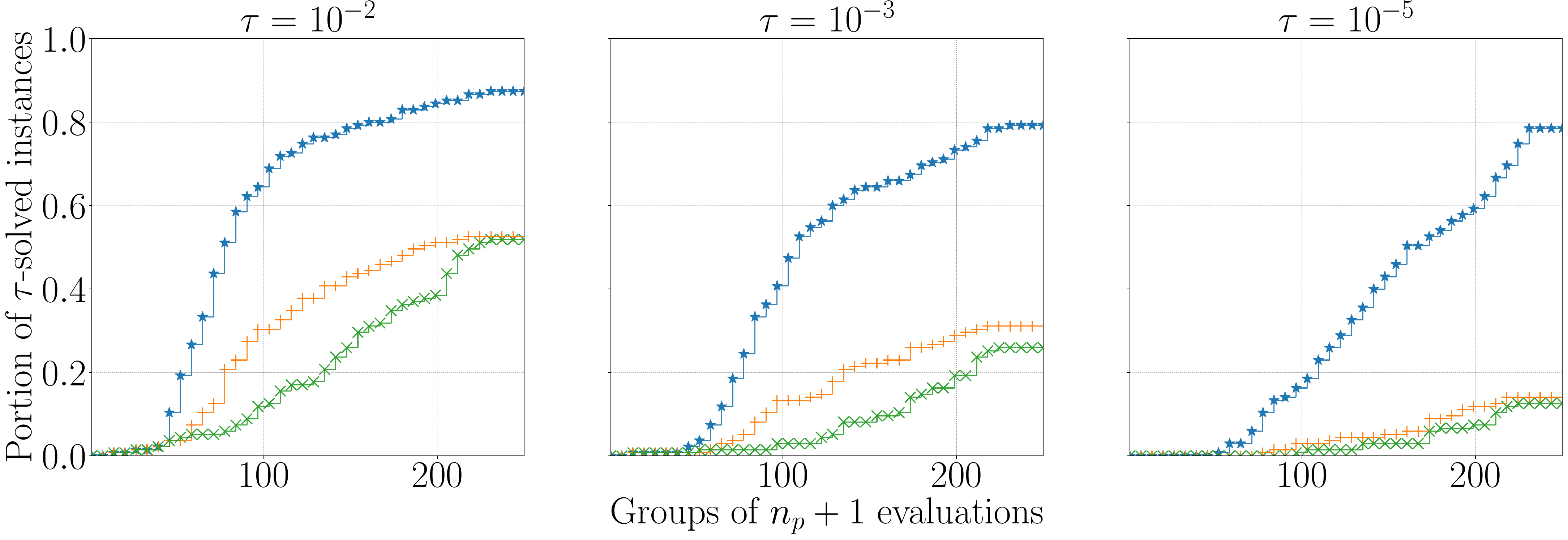}
\caption{Constrained test problems.}
\label{subfig:dataprofiles_constrained}
\end{subfigure}
\begin{subfigure}{\textwidth}
\centering
\includegraphics[width=0.55\linewidth]{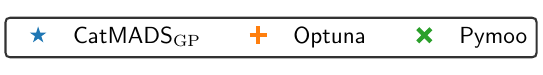}
\end{subfigure}
\caption{Comparison results on the test problems.}
\label{fig:data_profiles}
\end{figure}

The data profiles comparing \catmadsgp with the external solvers on the unconstrained and constrained instances are presented in \cref{subfig:dataprofiles_unconstrained,subfig:dataprofiles_constrained}, respectively.

For the unconstrained problems, \catmadsgp performs best for all three tolerances.
The difference increases for the stricter tolerances.
For $\tau=10^{-5}$, \catmadsgp solves approximately $70\%$ of the instances at $250(n_p+1)$, compared with about $20\%$ for {\sf Pymoo} and $10\%$ for {\sf Optuna}.

For the constrained problems, \catmadsgp also performs best for all three tolerances.
The difference is particularly pronounced for $\tau=10^{-3}$ and $\tau=10^{-5}$.
For instance, at $\tau=10^{-3}$ and $150(n_p+1)$, \catmadsgp solves approximately $60\%$ of the instances, compared with about $20\%$ for {\sf Optuna} and $10\%$ for {\sf Pymoo}.

%
%
%
%
%
%
%

\subsubsection{Mechanical-part application}
\label{sec:A}

{\blr
Data profiles provide a comparison over several mixed-variable problems. 
As shown in \cref{fig:histograms}, an important difficulty of the mechanical-part design problem is the identification of a feasible solution using a limited number of evaluations.
This problem is particularly relevant for assessing the performance of the proposed constrained surrogate-based neighborhoods.
}

For this focused study, the design problem is treated as a mixed-variable blackbox optimization problem.
To emulate a realistic blackbox setting with expensive evaluations, a small budget of $100$ evaluations is allowed.
The experiment is repeated over $50$ different seeds.
For each seed, all solvers are initialized with the same DoE comprising $20$ evaluations, corresponding to $20\%$ of the budget.
{\blr
The mean and standard deviation in \cref{tab:mechanical_part_results} are computed using only the runs in which a feasible solution is found.
The best-ranked runs correspond to the runs in which a solver obtains the lowest feasible objective value.
}

\begin{table}[htb!]
\centering
\normalsize
\begin{tabular}{lcccc}
\hline
\noalign{\vskip 1ex}
\raisebox{1.1ex}{Solver}
& \shortstack{Feasible \\ runs (\%)}
& \shortstack{Mean best feasible \\ objective}
& \shortstack{Standard \\ deviation}
& \shortstack{Best-ranked \\ runs (\%)} \\[1ex]
\hline
{\sf Optuna} & 82\% & 23.75 & 5.15 & 14\% \\[0.3ex]
{\sf Pymoo} & 54\% & 26.25 & 2.16 & 4\% \\[0.4ex]
\catmadsgp & \textbf{84\%} & \textbf{9.22} & 2.38 & \textbf{80\%} \\
\hline
\end{tabular}
\caption{Mechanical-part application results over 50 independent runs.}
\label{tab:mechanical_part_results}
\end{table}

{\blr
\cref{tab:mechanical_part_results} shows that \catmadsgp finds a feasible solution in $84\%$ of the runs, compared with $82\%$ for {\sf Optuna} and $54\%$ for {\sf Pymoo}.
Despite the similar feasibility rates of \catmadsgp and {\sf Optuna}, \catmadsgp obtains a much lower mean best feasible objective and is best-ranked in $80\%$ of the runs, compared with $14\%$ for {\sf Optuna} and $4\%$ for {\sf Pymoo}.
In the remaining $2\%$, none of the solvers finds a feasible solution.
}

%
\section{Discussion}
\label{sec:discussion}

{\blr This work introduces a novel approach for constructing categorical neighborhoods from surrogate models.}
These neighborhoods incorporate information from both the objective and the constraint functions, making them suitable for constrained mixed-variable optimization.
To balance similarity with respect to the objective and constraint functions, the construction relies on dominance relations.

As a proof of concept, the surrogate-based neighborhoods are first illustrated on the mechanical-part design problem.
A distinctive feature of this problem is the difficulty of identifying feasible points among the categorical components.
The proposed neighborhoods navigate this challenge more effectively than alternative neighborhood strategies.

The surrogate-based neighborhoods are then integrated into the \catmads algorithm, resulting in the new method \catmadsgp.
%
It is compared with two external mixed-variable solvers.
Data profiles indicate that \catmadsgp performs best on the considered \catsuite instances.
These results suggest that surrogate-based neighborhoods provide a promising mechanism for constrained mixed-variable blackbox optimization.

Further improvements could extend this study.
The current implementation of \catmadsgp uses a fixed number of neighbors throughout the iterations, and more sophisticated strategies could dynamically adjust the number of neighbors.
{\blr
The current neighborhood construction also requires enumerating the complete categorical domain, which may become impractical for very large categorical domains.
Constructing neighborhoods from only a partial enumeration of the categorical domain is therefore an interesting direction for future work.
}
Future extensions would define and analyze neighborhoods involving meta variables~\cite{G-2022-11}, \textit{i.e.}, variables whose values dictate the number of variables and/or constraints present in the optimization problem.



%



\section*{Data availability statement}
{\blr
The \catmadsgp implementation is publicly available through
\nomadbbo~\cite{NomadBBO}, currently distributed as a beta release at
\url{https://test.pypi.org/project/NomadBBO/}. 
The \catsuite benchmark problems are publicly available at
\url{https://github.com/bbopt/Cat-Suite}.
}

\section*{Conflict of interest statement}
The authors state that there are no conflicts of interest.


\section*{Acknowledgments} 
We gratefully acknowledge Dr. Paul Saves for his help with the implementation of SMT and the kernel. 

\bibliographystyle{plainurl}
\bibliography{bibliography}
\pdfbookmark[1]{References}{sec-refs}

\end{document}